\documentclass[11pt,a4paper,twoside]{article}
\usepackage{bcc20,subeqnarray}
\usepackage{bm}    
\usepackage{cite}  
\usepackage{epic,eepic}

\bcctitle{The multivariate Tutte polynomial (alias Potts model)
                for graphs and matroids}
\bccshorttitle{The multivariate Tutte polynomial}

\bccname{Alan D.~Sokal}
\bccshortname{A.D.~Sokal}

\bccaddressa{Department of Physics \\
New York University \\
4 Washington Place \\
New York, NY 10003  \\
USA}

\bccemaila{ {\tt sokal@nyu.edu} }

\newmathop{\deven}{{\rm DE}}
\newmathop{\dodd}{{\rm DO}}
\newmathop{\aut}{{\rm Aut}}
\newmathbin{\wreath}{{\rm wr}}

\newcommand{\be}{\begin{equation}}
\newcommand{\ee}{\end{equation}}
\newcommand{\<}{\langle}
\renewcommand{\>}{\rangle}

\def\reff#1{(\protect\ref{#1})}
\def\spose#1{\hbox to 0pt{#1\hss}}
\def\ltapprox{\mathrel{\spose{\lower 3pt\hbox{$\mathchar"218$}}
 \raise 2.0pt\hbox{$\mathchar"13C$}}}
\def\gtapprox{\mathrel{\spose{\lower 3pt\hbox{$\mathchar"218$}}
 \raise 2.0pt\hbox{$\mathchar"13E$}}}

\newcommand{\real}{\mathop{\rm Re}\nolimits}

\renewcommand{\diag}{\mathop{\rm diag}\nolimits}

\newcommand{\supp}{\mathop{\rm supp}\nolimits}
\newcommand{\per}{\mathop{\rm per}\nolimits}

\newcommand{\restrict}{\upharpoonright}
\newcommand{\drop}{\setminus}
\renewcommand{\emptyset}{\varnothing}

\newcommand{\scra}{{\mathcal{A}}}
\newcommand{\scrb}{{\mathcal{B}}}

\newcommand{\scrd}{{\mathcal{D}}}

\newcommand{\scrg}{{\mathcal{G}}}

\newcommand{\scrm}{{\mathcal{M}}}

\newcommand{\scrs}{{\mathcal{S}}}

\newcommand{\scrv}{{\mathcal{V}}}

\newcommand{\scry}{{\mathcal{Y}}}

\newcommand{\C}{{\mathbb C}}
\renewcommand{\Z}{{\mathbb Z}}

\newcommand{\Q}{{\mathbb Q}}

\newcommand{\Ztilde}{\widetilde{Z}}
\newcommand{\Ptilde}{\widetilde{P}}

\newcommand{\seriesq}{{\,\bowtie_q\,}}
\newcommand{\seriesnoq}{{\,\bowtie\,}}
\renewcommand{\parallel}{\Vert}

\newcommand{\bp}{ {\bf p} }
\newcommand{\bu}{ {\bf u} }
\newcommand{\bv}{ {\bf v} }
\newcommand{\bw}{ {\bf w} }
\newcommand{\bx}{ {\bf x} }
\def\qvbf{{q, \bv}}

\def\zgxconny{Z_G^{(x \leftrightarrow y)}}
\def\zgxnoconny{Z_G^{(x \not\leftrightarrow y)}}
\newcommand{\bvarphi}{ {\boldsymbol{\varphi}} }

\def\bsigma{{\boldsymbol{\sigma}}}

\def\bGamma{{\boldsymbol{\Gamma}}}

\def\psibar{{\bar{\psi}}}

\newenvironment{scarray}{
          \textfont0=\scriptfont0
          \scriptfont0=\scriptscriptfont0
          \textfont1=\scriptfont1
          \scriptfont1=\scriptscriptfont1
          \textfont2=\scriptfont2
          \scriptfont2=\scriptscriptfont2
          \textfont3=\scriptfont3
          \scriptfont3=\scriptscriptfont3
        
        \begin{array}{c}}{\end{array}}

\begin{document}
\makebcctitle

\begin{abstract}
The multivariate Tutte polynomial (known to physicists as the
Potts-model partition function) can be defined on an arbitrary
finite graph $G$, or more generally on an arbitrary matroid $M$,
and encodes much important combinatorial information about the graph
(indeed, in the matroid case it encodes the full structure of the
matroid).  It contains as a special case the familiar two-variable
Tutte polynomial --- and therefore also its one-variable specializations
such as the chromatic polynomial, the flow polynomial
and the reliability polynomial --- but is considerably more flexible.
I~begin by giving an introduction to all these problems,
stressing the advantages of working with the multivariate version.
I~then discuss some questions concerning the complex zeros
of the multivariate Tutte polynomial,
along with their physical interpretations in statistical mechanics
(in connection with the Yang--Lee approach to phase transitions)
and electrical circuit theory.
Along the way I~mention numerous open problems.
This survey is intended to be understandable to mathematicians with
no prior knowledge of physics.
\end{abstract}

\section{Introduction}  \label{sec1}

Let $G = (V,E)$ be a finite undirected graph
with vertex set $V$ and edge set $E$.\footnote{
   In this paper a ``graph'' is allowed to have
   loops and/or multiple edges unless explicitly stated otherwise.
}
The \defword{multivariate Tutte polynomial} of $G$ is, by definition,
the polynomial
\begin{equation}
   Z_G(q, \bv)   \;=\;
   \sum_{A \subseteq E}  q^{k(A)}  \prod_{e \in A}  v_e
   \;,
  \label{eq1.1}
\end{equation}
where $q$ and $\bv = \{v_e\}_{e \in E}$ are commuting indeterminates,
and $k(A)$ denotes the number of connected components in the subgraph $(V,A)$.
[It is sometimes convenient to consider instead
\begin{equation}
   \Ztilde_G(q, \bv)   \;\equiv\;
   q^{-|V|} Z_G(q, \bv)  \;=\;
   \sum_{A \subseteq E}  q^{k(A) - |V|}  \prod_{e \in A}  v_e
   \;,
  \label{eq.def.Ztilde}
\end{equation}
which is a polynomial in $q^{-1}$ and $\{v_e\}$.]
{}From a combinatorial point of view,
$Z_G$ is simply the multivariate generating polynomial
that enumerates the spanning subgraphs of $G$ according to
their precise edge content (with weight $v_e$ for the edge~$e$)
and their number of connected components (with weight $q$ for each component).
As we shall see, $Z_G$ encodes a vast amount of
combinatorial information about the graph $G$,
and contains many other well-known graph polynomials as special cases.
I~shall most often take an analytic point of view,
and treat $q$ and $\{v_e\}$ as real or complex variables.

In statistical physics, $Z_G$ is known as the
partition function of the \defword{$q$-state Potts model}\/
\cite{Potts_52,Wu_82,Wu_84},
for reasons to be explained in the next section.\footnote{
   The Potts model \cite{Potts_52} was invented in the early 1950s
   by Potts' thesis advisor Domb: see \cite{Domb_74}.
   The $q=2$ case, known as the Ising model \cite{Ising_25},
   was invented in 1920 by Ising's thesis advisor Lenz \cite{Lenz_20}:
   see \cite{Brush_67,Kobe_97,Stutz_99} for a fascinating history.
   (I~hasten to add that these are the only two cases I~know of
   where the thesis advisor's invention was named after the graduate student,
   rather than the other way around.)
   The $q=4$ case, which is a special case of the Ashkin--Teller model,
   was invented in 1943 by Ashkin and Teller \cite{Ashkin-Teller_43}.
   For a review of the Potts model, see \cite{Wu_82,Wu_84}.
}
The Potts model --- along with its close relative,
the Fortuin--Kasteleyn random-cluster model
\cite{Kasteleyn_69,Fortuin_72,Grimmett_book_to_appear} ---
plays an important role in the theory of phase transitions
and critical phenomena \cite{Baxter_82,Nienhuis_84,DiFrancesco_97},
and physicists have developed a significant amount
of information (both rigorous and non-rigorous) about its properties.
As will be explained in Section~\ref{sec.complex},
one way of analyzing phase transitions is to study
the zeros of $Z_G(q, \bv)$ when $q$ and/or $\{v_e\}$
are treated as {\em complex}\/ variables.

If we set all the edge weights $v_e$ equal to the same value $v$,
we obtain a two-variable polynomial $Z_G(q,v)$ that is essentially
equivalent to the standard Tutte polynomial $T_G(x,y)$
[see Section~\ref{subsec.tutte} for details].
But the main message of this paper is that
it is often useful to consider the multivariate polynomial
$Z_G(q, \bv)$, even if one is ultimately interested in a particular
two-variable or one-variable specialization.
For instance, $Z_G(q, \bv)$ is \defword{multiaffine}
in the variables $\bv$ (i.e., of degree 1 in each $v_e$ separately);
often a multiaffine polynomial in many variables is easier to handle
than a general polynomial in a single variable
(e.g., it may permit simple proofs by induction on the number of variables).
Furthermore, many natural operations on graphs,
such as the reduction of edges in series or parallel,
lead out of the class of ``all $v_e$ equal''.
I~shall illustrate the advantages of this ``multivariate ideology''
by showing several instances in which the multivariate extension
of a single-variable result is not only vastly more powerful
but also much \emph{easier} to prove:
indeed, in one case, a 20-page proof is reduced to a few lines.

All of these considerations can be extended from graphs
to matroids\footnote{
   See Oxley \cite{Oxley_92} for an excellent introduction to
   matroid theory.
}:
if $M$ is a matroid with ground set $E$ and rank function $r_M$,
its multivariate Tutte polynomial is defined by
\begin{equation}
   \Ztilde_M(q, \bv)   \;=\;
   \sum_{A \subseteq E}  q^{-r_M(A)}  \prod_{e \in A}  v_e
   \;,
  \label{eq.defZ.matroid}
\end{equation}
which is a polynomial in $q^{-1}$ and $\{v_e\}$.
This extends the graph definition in the sense that if $G$ is a graph
and $M(G)$ is its cycle matroid, then
\begin{equation}
   \Ztilde_{M(G)}(q, \bv)   \;=\;
   \Ztilde_G(q, \bv)
\end{equation}
[because $r_{M(G)}(A) = |V| - k(A)$].
Since a matroid is completely determined by its rank function,
$\Ztilde_M$ is simply an algebraic encoding of {\em all}\/
the information about the matroid $M$.
Moreover, my earlier statement that $Z_G$ encodes ``a vast amount''
of information about the graph $G$ can now be made more precise:
$Z_G$ encodes the number of vertices $|V|$
together with all the information about $G$
that is contained in its cycle matroid $M(G)$ [and no other information].

I~am now convinced that matroids are the natural category
for studying the multivariate Tutte polynomial.
Of course, results that hold for graphic matroids
may or may not hold for larger classes of matroids.

\section{The multivariate Tutte polynomial for graphs, and its specializations}
   \label{sec2}

First, a trivial but useful observation:
Starting from the definition
\begin{equation}
   Z_G(q, \bv)   \;=\;
   \sum_{A \subseteq E}  q^{k(A)}  \prod_{e \in A}  v_e
 \label{eq2.1}
\end{equation}
and using the simple relation
\begin{equation}
   k(A)   \;=\;   |V| \,-\, |A| \,+\, c(A)
 \label{eq.cyclomatic}
\end{equation}
where $c(A)$ is the cyclomatic number
(i.e., number of linearly independent cycles) of the graph $(V,A)$,
we can rewrite $Z_G$ in the alternate form
\begin{equation}
   Z_G(q, \bv)   \;=\;
   q^{|V|} \, \sum_{ A \subseteq E }  q^{c(A)} \prod_{e \in A} {v_e \over q}
   \;.
 \label{eq1.1cycles}
\end{equation}

\begin{example}
 \label{sec.potts.defs.trees}
{\em Trees.}\/
For any tree $T=(V,E)$ we have
\begin{equation}
   Z_T(q, \bv)   \;=\;  q \, \prod_{e \in E} (q+v_e)
   \;.
\end{equation}
This follows immediately from \reff{eq1.1cycles},
using the fact that $c(A) = 0$ for all $A$.
$\quad\square$
\end{example}

\begin{example}
 \label{sec.potts.defs.cycles}
{\em Cycles.}\/
Let $G=(V,E)$ be the cycle $C_n$.  Then
\begin{equation}
   Z_{C_n}(q, \bv)   \;=\;
       \prod_{e \in E}  (q+v_e)   \,+\, (q-1) \prod_{e \in E}  v_e
   \;.
\end{equation}
This follows immediately from \reff{eq1.1cycles},
using the fact that $c(A) = 0$ for all $A \subsetneq E$
and $c(E) = 1$.
$\quad\square$
\end{example}

Let us now consider some of the polynomials that can be obtained
from $Z_G(q, \bv)$ by specializing the value of $q$:

\subsection{\mbox{\protect\boldmath $q=1$}}

When $q=1$, the multivariate Tutte polynomial becomes trivial:
\begin{equation}
   Z_G(1, \bv)  \;=\;  \prod_{e \in E} (1+v_e)
   \;.
\end{equation}

\subsection{\mbox{\protect\boldmath $q=1,2,3,\ldots$ ($q$-state Potts model)
   and the chromatic polynomial}}
  \label{subsec.FK.graphs}

Let $q$ be a positive integer;
then the \defword{$q$-state Potts-model partition function}
for the graph $G$ is defined by
\begin{equation}
   Z_G^{\rm Potts}(\qvbf)   \;=\;
   \sum_{ \sigma \colon\, V \to \{ 1,2,\ldots,q \} }
   \; \prod_{e \in E}  \,
      \biggl[ 1 + v_e \delta(\sigma_{x_1(e)}, \sigma_{x_2(e)}) \biggr]
   \;.
 \label{def.ZPotts}
\end{equation}
Here the sum runs over all maps $\sigma\colon\, V \to \{ 1,2,\ldots,q \}$,
and we sometimes write $\sigma_x$ as a synonym for $\sigma(x)$;
the $\delta$ is the Kronecker delta
\begin{equation}
   \delta(a,b)   \;=\;    \begin{cases}
                              1  & \text{if $a=b$} \\
                              0  & \text{if $a \neq b$}
                          \end{cases}
\end{equation}
and $x_1(e), x_2(e) \in V$ are the two endpoints of the edge $e$
(in arbitrary order).
We usually consider $v_e$ to be a real or complex variable.

In statistical physics, this formula arises as follows:
In the Potts model \cite{Potts_52,Wu_82,Wu_84},
an ``atom'' (or ``spin'') at the site $x \in V$ can exist in any one of
$q$ different states.
The {\em energy}\/ of a configuration is the sum, over all edges $e \in E$,
of $0$ if the spins at the two endpoints of that edge are unequal
and $-J_e$ if they are equal.
The {\em Boltzmann weight}\/ of a configuration is then $e^{-\beta H}$,
where $H$ is the energy of the configuration
and $\beta \ge 0$ is the inverse temperature.
The {\em partition function}\/ is the sum, over all configurations,
of their Boltzmann weights.
Clearly this is just a rephrasing of \reff{def.ZPotts},
with $v_e = e^{\beta J_e} - 1$.
A coupling $J_e$ (or $v_e$)
is called {\em ferromagnetic}\/ if $J_e \ge 0$ ($v_e \ge 0$),
as it is then favored for adjacent spins to take the same value;
{\em antiferromagnetic}\/ if $-\infty \le J_e \le 0$ ($-1 \le v_e \le 0$),
as it is then favored for adjacent spins to take different values;
and {\em unphysical}\/ if $v_e \notin [-1,\infty)$,
as the weights are then no longer nonnegative.

It is far from obvious that $Z_G^{\rm Potts}(\qvbf)$,
which is defined separately for each positive integer $q$,
is in fact the restriction to $q \in \Z_+$
of a {\em polynomial}\/ in $q$.
But this is in fact the case, and indeed we have:

\begin{theorem}[Fortuin--Kasteleyn representation of the Potts model]
   \label{thm.FK}
\hfill\break
\vspace*{-0.4cm}
\par\noindent
For integer $q \ge 1$,
\begin{equation}
   Z_G^{\rm Potts}(\qvbf) \;=\;  Z_G(q, \bv)  \;.
 \label{eq.FK.identity}
\end{equation}
That is, the Potts-model partition function
is simply the specialization of the multivariate Tutte polynomial
to $q \in \Z_+$.
\end{theorem}

\begin{proof}
In \reff{def.ZPotts}, expand out the product over $e \in E$,
and let $A \subseteq E$ be the set of edges for which the term
$v_e \delta(\sigma_{x_1(e)}, \sigma_{x_2(e)})$ is taken.
Now perform the sum over configurations $\{ \sigma_x \}_{x \in V}$:
in each component of the subgraph $(V,A)$
the color $\sigma_x$ must be constant,
and there are no other constraints.
Therefore,
\begin{equation}
   Z_G(\qvbf) \;=\;
   \sum_{ A \subseteq E }  q^{k(A)}  \prod_{e \in A}  v_e
   \;,
  \label{eq1.1bis}
\end{equation}
as was to be proved.
\end{proof}

The subgraph expansion \reff{eq1.1bis}
of the Potts-model partition function was discovered in the late 1960s by
Fortuin and Kasteleyn \cite{Kasteleyn_69,Fortuin_72}.\footnote{
   In the special case $v_e = -1$ it was discovered many decades earlier by
   Birkhoff \cite{Birkhoff_12} and Whitney \cite{Whitney_32a};
   see also Tutte \cite{Tutte_47,Tutte_54}.
}
Please note that if $q \ge 0$ and $v_e \ge 0$ for all $e$,
then the weights in \reff{eq1.1bis} are nonnegative,
and so can be interpreted probabilistically (after normalization by $Z_G$).
The resulting probability measure on $2^E$ is called the
\defword{Fortuin--Kasteleyn random-cluster model}\/
\cite{Kasteleyn_69,Fortuin_72,Grimmett_book_to_appear}
for the graph $G$.
Furthermore, for $q \in \Z_+$ a natural joint probability measure
on spin states and edge occupations
[i.e., a coupling of \reff{def.ZPotts} and \reff{eq1.1bis}]
has been given by Edwards and Sokal \cite{Edwards-Sokal};
it underlies the Swendsen--Wang \cite{Swendsen_87} Monte Carlo algorithm
for the ferromagnetic Potts model,
and is also useful in the rigorous analysis of the Potts model.

Another important special case arises when $v_e = -1$ for all edges $e$:
then $Z_G^{\rm Potts}$ gives weight 1 to each proper coloring
and weight 0 to each improper coloring,
and so counts the proper colorings.
It follows from Theorem~\ref{thm.FK} that the number of proper $q$-colorings
of $G$ is in fact the restriction to $q \in \Z_+$
of a polynomial in $q$, called the \defword{chromatic polynomial}
$P_G(q) = Z_G(q,-1)$.
The chromatic polynomial thus corresponds to
the zero-temperature ($\beta \to +\infty$) limit
of the antiferromagnetic ($J_e < 0$) Potts model.
Many properties of the chromatic polynomial in fact extend to the
entire antiferromagnetic region (i.e., $-1 \le v_e \le 0$ for all $e$).

\subsection{\mbox{\protect\boldmath $q \to 0$} limits}   \label{subsec.q0}

Let us now consider the different ways in which a meaningful
$q \to 0$ limit can be taken in the $q$-state Potts model.

The simplest limit is to take $q \to 0$ with fixed $\bv$.
{}From the definition \reff{eq2.1} we see that
this selects out the subgraphs $A \subseteq E$
having the smallest possible number of connected components;
the minimum achievable value is of course $k(G)$ itself
(= 1 in case $G$ is connected).
We therefore have
\begin{equation}
   \lim_{q \to 0} q^{-k(G)} Z_G(q,\bv) \;=\;  C_G(\bv)
   \;,
\end{equation}
where
\begin{equation}
   C_G(\bv) \;=\;  \sum\limits_{\begin{scarray}
                                       A \subseteq E \\
                                       k(A) = k(G)
                                    \end{scarray}}
                       \prod_{e \in A}  v_e
\end{equation}
is the generating polynomial of ``maximally connected spanning subgraphs''
(= \defword{connected spanning subgraphs}\/ in case $G$ is connected).

A different limit can be obtained by taking $q \to 0$
with fixed values of $\bw = \bv/q$.
{}From \reff{eq1.1cycles} we see that this selects out
the subgraphs $A \subseteq E$ having the smallest possible cyclomatic number;
the minimum achievable value is of course 0.
We therefore have \cite{Stephen_76,Wu_77}
\begin{equation}
   \lim_{q \to 0} q^{-|V|} Z_G(q,q\bw) \;=\;  F_G(\bw)
   \;,
 \label{eq.limit.FG}
\end{equation}
where
\begin{equation}
   F_G(\bw) \;=\;  \sum\limits_{\begin{scarray}
                                       A \subseteq E \\
                                       c(A) = 0
                                    \end{scarray}}
                       \prod_{e \in A}  w_e
   \label{def_F}
\end{equation}
is the generating polynomial of {\em spanning forests}\/.

Finally, suppose that in $C_G(\bv)$
we replace each edge weight $v_e$ by $\lambda v_e$
and then take $\lambda \to 0$.
This obviously selects out,
from among the maximally connected spanning subgraphs,
those having the fewest edges:
these are precisely the maximal spanning forests
(= \defword{spanning trees}\/ in case $G$ is connected),
and they all have exactly $|V| - k(G)$ edges.
Hence
\begin{equation}
   \lim_{\lambda \to 0} \lambda^{k(G)-|V|} C_G(\lambda \bv) \;=\;
   T_G(\bv)
   \;,
\end{equation}
where
\begin{equation}
   T_G(\bv) \;=\;  \sum\limits_{\begin{scarray}
                                       A \subseteq E \\
                                       k(A) = k(G) \\
                                       c(A) = 0
                                    \end{scarray}}
                       \prod_{e \in A}  v_e
 \label{def.tree}
\end{equation}
is the generating polynomial of maximal spanning forests.\footnote{
   I trust that there will be no confusion between the
   generating polynomial $T_G(\bv)$
   and the Tutte polynomial $T_G(x,y)$.
   I have used here the letter $T$ because in the most important applications
   the graph $G$ is connected, so that $T_G(\bv)$ is the generating
   polynomial of spanning {\em trees}\/.
}
Alternatively, suppose that in $F_G(\bw)$
we replace each edge weight $w_e$ by $\lambda w_e$
and then take $\lambda \to \infty$.
This obviously selects out, from among the spanning forests,
those having the greatest number of edges:
these are once again the maximal spanning forests.
Hence
\begin{equation}
   \lim_{\lambda \to \infty} \lambda^{k(G)-|V|} F_G(\lambda \bw) \;=\;
   T_G(\bw)
   \;.
\end{equation}

In summary, we have the following scheme for the $q \to 0$ limits
of the multivariate Tutte polynomial:
\begin{equation}
  \begin{array}{ccccc}
        &   &  C_G(\bv)   &     &    \\  
        &   \setlength{\unitlength}{1pt}
            \begin{picture}(70,50)(0,0)
            \drawline(0,0)(25.2,18)
            \put(0,23){{\small $q \to 0,\, \bv \, {\rm fixed}$}}
            \drawline(49,35)(70,50)
            \drawline(70,50)(64.34,41.76)
            \drawline(70,50)(60.37,47.32)
            \end{picture} &
        &   \setlength{\unitlength}{1pt}
            \begin{picture}(70,50)(0,0)
            \drawline(0,50)(21,35)
            \put(20,23){{\small $\bv \, {\rm infinitesimal}$}}
            \drawline(44.8,18)(70,0)
            \drawline(70,0)(60.37,2.68)
            \drawline(70,0)(64.34,8.24)
            \end{picture} &              \\
   Z_G(q,\bv)  &    &      &    &   T_G(\bv \: {\rm or} \: \bw) \\
        &   \setlength{\unitlength}{1pt}
            \begin{picture}(70,50)(0,0)
            \drawline(0,50)(21,35)
            \put(0,23){{\small $q \to 0,\, \bw = \bv/q \, {\rm fixed}$}}
            \drawline(44.8,18)(70,0)
            \drawline(70,0)(60.37,2.68)
            \drawline(70,0)(64.34,8.24)
            \end{picture} &
        &   \setlength{\unitlength}{1pt}
            \begin{picture}(70,50)(0,0)
            \drawline(0,0)(25.2,18)
            \put(20,23){{\small $\bw \, {\rm infinite}$}}
            \drawline(49,35)(70,50)
            \drawline(70,50)(64.34,41.76)
            \drawline(70,50)(60.37,47.32)
            \end{picture} & \\
        &   &  F_G(\bw)   &     &
  \end{array}
 \label{scheme.qto0}
\end{equation}

Finally, maximal spanning forests (= spanning trees in case $G$ is connected)
can also be obtained directly from $Z_G(q,\bv)$ by a one-step process
in which the limit $q \to 0$ is taken at fixed $\bx = \bv/q^\alpha$,
where $0 < \alpha < 1$ \cite{Fortuin_72,Stephen_76,Wu_77,Haggstrom_95}.
Indeed, simple manipulation of \reff{eq2.1} and \reff{eq.cyclomatic} yields
\begin{equation}
   Z_G(q, q^\alpha \bx) \;=\;
   q^{\alpha|V|} \sum_{ A \subseteq E }  q^{\alpha c(A) + (1-\alpha)k(A)}
       \prod_{e \in A} x_e
   \;.
 \label{eq1.1alpha}
\end{equation}
The quantity $\alpha c(A) + (1-\alpha)k(A)$ is minimized on (and only on)
maximal spanning forests, where it takes the value $(1-\alpha)k(G)$.
Hence
\begin{equation}
   \lim_{q \to 0} q^{-\alpha|V| - (1-\alpha)k(G)} Z_G(q, q^\alpha \bx)
   \;=\;  T_G(\bx)
   \;.
\end{equation}


\subsection{The multivariate flow polynomial}  \label{subsec.flow}

Let $G = (V,E)$ be a finite undirected graph,
and let $\Gamma$ be a finite abelian group of order $q = |\Gamma|$.
Let us choose arbitrarily an orientation for each edge of $G$
(all subsequent results will be independent of this choice).
A \defword{$\Gamma$-flow} on $G$
is a mapping $\psi\colon\, E \to \Gamma$ that satisfies current conservation
at every vertex.
A $\Gamma$-flow $\psi$ on $G$ is said to be \defword{nowhere-zero}
in case $\psi(e) \neq 0$ for all $e \in E$.

Let $F_G(\Gamma)$ be the number of
nowhere-zero $\Gamma$-flows on $G$.\footnote{
   I trust that there will be no confusion between the
   flow polynomial $F_G(\Gamma)$ [and its relatives to be defined below]
   and the spanning-forest polynomial $F_G(\bw)$.
}
It is a well-known (but at first sight quite surprising) fact
that $F_G(\Gamma)$ depends only on the order of $\Gamma$
and not on any other aspect of the structure of $\Gamma$;
so we can write it as $F_G(q)$.
Furthermore, it turns out that $F_G(q)$ is the restriction
to $\Z_+$ of a polynomial in $q$,
called the \defword{flow polynomial}\/ of $G$
and also written $F_G(q)$.


In order to prove these facts (and much more),
let us generalize $F_G(\Gamma)$ as follows:
Assign to each edge $e \in E$ a weight $u_e$,
and write $\bu = \{u_e\}_{e \in E}$.
Then define
\begin{equation}
   F_G(\Gamma, \bu)   \;=\;
   \sum_{\Gamma\hbox{\scriptsize -flows}\, \psi}  \;\,  \prod_{e \in E}  \,
      \biggl[ 1 + u_e \delta(\psi(e),0) \biggr]
   \;,
 \label{eq.sec.flowpoly.1}
\end{equation}
where the sum runs over all $\Gamma$-flows $\psi$ on $G$.
If we take $u_e = -1$ for all $e$, this reduces to our friend $F_G(\Gamma)$.
We now assert that $F_G(\Gamma, \bu)$ depends only
on $q = |\Gamma|$ --- so we can write it as $F_G(q, \bu)$ ---
and that $F_G(q, \bu)$ is in fact the restriction to $q \in \Z_+$ of
a polynomial in $q$ and $\{u_e\}$,
which we shall call the \defword{multivariate flow polynomial}\/.
Indeed, the multivariate flow polynomial is nothing other than
the multivariate Tutte polynomial in disguise:

\begin{theorem}
   \label{thm.potts.flowpoly}
Let $G=(V,E)$ be a finite undirected graph,
and let $\Gamma$ be a finite abelian group of order $q = |\Gamma|$.
Then
\begin{equation}
   F_G(\Gamma, \bu)   \;=\;
   q^{-|V|} \, \left( \prod\limits_{e \in E} u_e \right) \, Z_G(q, q/\bu)
   \label{thm.potts.flowpoly.eq1}
\end{equation}
(where division is of course understood edge-wise)
and in particular
\begin{equation}
    F_G(\Gamma)   \;=\;  q^{-|V|} \, (-1)^{|E|} \, Z_G(q,-q)
    \;.
\end{equation}
Consequently, $F_G(\Gamma)$ and $F_G(\Gamma, \bu)$ depend only on
$q = |\Gamma|$ and not on any other aspect of the structure of $\Gamma$.
\end{theorem}

\begin{proof}
Note first that the total number of $\Gamma$-flows
on a graph $G$ is $|\Gamma|^{c(G)}$,
where $c(G)$ is the cyclomatic number of $G$.
Moreover, the $\Gamma$-flows that vanish on a subset $A \subseteq E$
are in one-to-one correspondence with the $\Gamma$-flows on the
subgraph $(V, E \setminus A)$,
so their number is $|\Gamma|^{c(E \setminus A)}$.
Expanding the product in \reff{eq.sec.flowpoly.1}, we find
\begin{subeqnarray}
   F_G(\Gamma, \bu)   & = &
      \sum\limits_{A \subseteq E} q^{c(E \setminus A)}  \prod_{e \in A} u_e
     \\[2mm]
   & = &
   \left( \prod\limits_{e \in E} u_e \right) \,
      \sum\limits_{A' \subseteq E} q^{c(A')}  \prod_{e \in A'} {1 \over u_e}
     \\[2mm]
   & = &
   q^{-|V|} \, \left( \prod\limits_{e \in E} u_e \right) \,
      Z_G(q, q/\bu)
\end{subeqnarray}
by comparison with \reff{eq1.1cycles}.
\end{proof}

Theorem~\ref{thm.potts.flowpoly} is in fact a special case
of a more general identity
\cite{Wu-Wang_76,Savit_80,Druhl_82,Caracciolo-Sportiello_04}
relating ``spin models'' taking values in an abelian group $\Gamma$
to ``flow models'' taking values in the dual group $\Gamma^*$:
here the Boltzmann weight $W_e(\sigma_x - \sigma_y)$ for an edge $e=xy$
in the spin model maps onto
a weight $\widehat{W}_e(\psi(e))$ in the flow model,
where $\widehat{\hphantom{W}}$ denotes Fourier transformation.
In the special case of the Potts model,
we have $W_e = {\bf 1} + v_e \delta$,
so that $\widehat{W}_e = v_e {\bf 1} + q \delta$.
Setting $u_e = q/v_e$ and relabelling $\Gamma^*$ as $\Gamma$,
we obtain \reff{thm.potts.flowpoly.eq1}.

\subsection{Comparison to the standard Tutte polynomial}
   \label{subsec.tutte}

The Tutte polynomial $T_G(x,y)$ is conventionally defined by
   \cite[p.~45]{Welsh_93} \cite[pp.~124--127]{Brylawski_92}
\begin{subeqnarray}
   T_G(x,y)   & = &
   \sum\limits_{A \subseteq E}  (x-1)^{r(E)-r(A)} \, (y-1)^{|A| - r(A)}
       \\[2mm]
   & = &
   \sum\limits_{A \subseteq E}  (x-1)^{k(A)-k(E)} \, (y-1)^{|A| + k(A) - |V|}
 \label{eq.Tutte.1}
\end{subeqnarray}
where $r(A) = |V| - k(A)$ is the rank of the set $A$
in the cycle matroid $M(G)$.
[Note also that $|A| - r(A)$ is the cyclomatic number $c(A)$.]
Comparison with \reff{eq1.1} yields
\begin{equation}
   T_G(x,y)   \;=\;   (x-1)^{-k(E)} \, (y-1)^{-|V|} \,
                         Z_G \bigl( (x-1)(y-1), \, y-1 \bigr)   \;.
 \label{eq.Tutte.2}
\end{equation}
In other words, $T_G(x,y)$ and $Z_G(q,v)$ are essentially equivalent
under the change of variables
\begin{subeqnarray}
   x  & = &  1 + q/v   \\
   y  & = &  1 + v     \\[4mm]
   q  & = &  (x-1)(y-1) \\
   v  & = &  y-1
 \label{eq.Tutte.3}
\end{subeqnarray}

One advantage of the Tutte notation is that it allows a slightly smoother
treatment of the $q \to 0$ limit,
in which the limiting processes employed in Section~\ref{subsec.q0}
can be replaced by simple evaluation at $x=1$ or $y=1$.
Thus, the univariate maximally-connected-spanning-subgraph
and spanning-forest polynomials are given, respectively, by
\begin{eqnarray}
   C_G(v)  & = &  v^{|V|-k(E)} \, T_G(1,1+v)  \\[2mm]
   F_G(w)  & = &  w^{|V|-k(E)} \, T_G(1+1/w,1)
\end{eqnarray}
A second advantage of the Tutte notation is that duality
[cf.\ \reff{eq.duality} below] takes the simple form
\begin{equation}
   T_{G^*}(x,y)  \;=\;  T_G(y,x)   \;.
 \label{eq.duality.Tutte}
\end{equation}

But the Tutte notation also has a severe disadvantage:
the use of the variables $x$ and $y$
conceals the fact that the particular combinations $q$ and $v$
play {\em very different roles}\/;
$q$ is a global variable,
while $v$ can be given separate values $v_e$ on each edge.
In particular, the Tutte notation {\em makes it impossible}\/
to assign unequal weights $v_e$ to the edges.
I~therefore recommend use of the notation $Z_G(q,\bv)$
whenever the ``multivariate ideology'' is potentially of use.

Let me conclude by observing that numerous specific evaluations of 
the Tutte polynomial have been given combinatorial interpretations,
as counting some set of objects associated to the graph $G$
(see e.g.\ \cite{Brylawski_92,Welsh_93}).
It would be an interesting project to seek to extend these
counting problems to ``counting with weights'',
i.e., to obtain suitably defined
univariate or multivariate generating polynomials
for the objects in question
as specializations of $Z_G(q,v)$ or $Z_G(q,\bv)$, respectively.
A few examples of this are given in
Sections~\ref{subsec.q0} and \ref{subsec.flow} above.

\section{The multivariate Tutte polynomial for matroids}

Many (but not all) of the properties of the multivariate Tutte polynomial
for graphs can be carried over to the matroid version \reff{eq.defZ.matroid}.
Where the treatment is essentially identical, I shall be brief.

\subsection{\mbox{\protect\boldmath $q=1$}}

When $q=1$, the multivariate Tutte polynomial becomes trivial:
\begin{equation}
   \Ztilde_M(1, \bv)  \;=\;  \prod_{e \in E} (1+v_e)
   \;.
\end{equation}

\subsection{\mbox{\protect\boldmath $q=1,2,3,\ldots$
              (generalized Potts models)}}

For {\em graphs}\/, we saw in Section~\ref{subsec.FK.graphs} that
the multivariate Tutte polynomial $Z_G(q,\bv)$ for integer $q \ge 1$
counts $q$-colorings of the vertices of $G$,
with weights that depend on whether each edge $e$ is
properly or improperly colored [cf.\ \reff{def.ZPotts}/\reff{eq.FK.identity}].
In particular, the chromatic polynomial $P_G(q) = Z_G(q,-1)$
counts proper $q$-colorings of $G$.
One consequence of this is that the nonnegative integer zeros of $P_G(q)$
are consecutive integers, namely $0,1,2,\ldots,\chi(G)-1$,
where $\chi(G)$ is the chromatic number of $G$.\footnote{
   In fact, these are {\em all}\/ the integer zeros of $P_G(q)$,
   as it is easy to show that $P_G(q)$ has no negative real zeros
   (when $G$ is loopless): see e.g.\ \cite{Read_88}.
}

Simple examples show that the latter property does not hold in general
for the chromatic polynomial
(also known as \defword{characteristic polynomial}\/)
of a matroid, $\Ptilde_M(q) = \Ztilde_M(q,-1)$.\footnote{
   The chromatic (or characteristic) polynomial of a matroid $M$
   is ordinarily defined as $P_M(q) = q^{r(M)} \Ptilde_M(q)$,
   so that it is a polynomial in nonnegative powers of $q$.
   This is a matter of taste.
}
For instance, for the uniform matroid $U_{2,n}$ we have
\begin{equation}
   \Ztilde_{U_{2,n}}(q,\bv)  \;=\;
       q^{-2} \prod_{i=1}^n (1+v_i)  \,+\, (1 - q^{-2}) \,+\,
                                 (q^{-1} - q^{-2}) \sum_{i=1}^n v_i
\end{equation}
and hence
\begin{equation}
   \Ptilde_{U_{2,n}}(q)  \;=\;  q^{-2} (q-1)(q+1-n)
   \;,
\end{equation}
whose roots are $q=1$ and $q=n-1$;
so the consecutive-root property is violated for all $n \ge 4$.
It therefore seems unlikely that $\Ztilde_M(q, \bv)$
can be given a ``coloring'' interpretation
for {\em arbitrary}\/ matroids $M$ and {\em arbitrary}\/ integers $q \ge 1$.
Nevertheless, in some special cases such a representation can be given,
as I would now like to sketch;
this account is drawn from work in progress with
Sergio Caracciolo and Andrea Sportiello \cite{gen_potts}.

The appropriate general context for this discussion
is that of \defword{abelian-group-valued statistical-mechanics models}\/,
defined as follows \cite{Caracciolo-Sportiello_04}:
Let $V$ and $E$ be finite index sets
(for \defword{variables}\/ and \defword{interactions}\/, respectively);
let $R$ be a ring with identity and let $\scrm$ be a unitary right $R$-module;
let $B = (b_{ie})_{i \in V,\, e \in E}$ be a matrix with elements in $R$;
let $\mu$ be Haar measure on $\scrm$, considered as an additive abelian group;
and for each $e \in E$, let $W_e \colon\, \scrm \to \C$ be a function.
We then define the partition function
\begin{equation}
   \Ztilde_{R,\scrm,B}({\bf W})  \;=\;
   \int \prod_{e \in E} W_e\biggl(\, \sum_{i \in V} \sigma_i b_{ie} \biggr)
        \;  \prod_{i \in V} d\mu(\sigma_i)
\end{equation}
where ${\bf W} = \{W_e\}_{e \in E}$.
(If $\scrm$ is infinite, then the collection of functions ${\bf W}$
 needs to satisfy suitable integrability conditions.)

[We can alternatively let $\scrm$ be an arbitrary
additive abelian group $\scrg$,
and let $R$ be some subring of
the ring ${\rm End}(\scrg)$ of endomorphisms of $\scrg$.
It~is easy to see that the two formulations are equivalent.]

Two important special cases are:
\begin{itemize}
   \item[1)]  $R = \Z$ and $\scrm$ is an arbitrary
       additive abelian group $\scrg$
       (on which $\Z$ acts in the obvious way).
   \item[2)]  $R$ is a field $F$, and $\scrm$ is a finite-dimensional
       vector space over $F$ (which can be taken, without loss of
       generality, to be $F^k$ for some integer $k \ge 0$).
\end{itemize}

A \defword{generalized Potts model}\/ is the special case of this setup
in which $\scrm$ is finite,
$\mu$ is normalized counting measure on $\scrm$
(i.e., each element of $\scrm$ gets weight $1/|\scrm|$),
and each $W_e$ is of the form
\begin{equation}
   W_e(\sigma)  \;=\;  1 \,+\, v_e \delta(\sigma,0)
\end{equation}
for some weights $\bv = \{v_e\}_{e \in E}$.  We then have
\begin{equation}
   \Ztilde^{\rm Potts}_{R,\scrm,B}(\bv)
   \;=\;
   |\scrm|^{-|V|} \! \sum_{ \sigma \colon V \to \scrm }
       \;
       \prod_{e \in E}
       \Biggl[ 1 + v_e \delta\biggl(\,\sum_{i \in V} \sigma_i b_{ie},\, 0
                             \biggr)
       \Biggr]
   \;.
 \label{def.gen.Potts}
\end{equation}
Note that if $G=(V,E)$ is a graph,
$B$ is the directed vertex-edge incidence matrix for some orientation of $G$
(hence $R = \Z$), and $\scrm$ is a finite abelian group
of order $q = |\scrm|$, we have
$\Ztilde^{\rm Potts}_{R,\scrm,B}(\bv) = \Ztilde^{\rm Potts}_G(q,\bv)
 \equiv q^{-|V|} Z^{\rm Potts}_G(q,\bv)$.
So the models \reff{def.gen.Potts} do indeed constitute a generalization
of the Potts models \reff{def.ZPotts}.
 
As far as I know, there does not exist,
{\em at this level of generality}\/,
a satisfactory analogue of the Fortuin--Kasteleyn representation
\reff{eq.FK.identity}.  But at least in some special cases
such a representation does exist \cite{gen_potts}:

\begin{theorem}[Fortuin--Kasteleyn representation for generalized Potts models]
   \label{thm.gen.FK}
\hfill\break
\vspace*{-0.8cm}
\begin{itemize}
   \item[(a)]  Let $M$ be a regular matroid, let $B$ be a totally unimodular
      matrix of integers representing $M$ over $\Q$,
      and let $\scrg$ be a finite abelian group of order $q = |\scrg|$.
      Then
\begin{equation}
     \Ztilde^{\rm Potts}_{\Z,\scrg,B}(\bv)  \;=\; \Ztilde_M(q,\bv)
     \;.
 \label{eq.gen.FK.identity1}
\end{equation}
   \item[(b)]  Let $M$ be a matroid, let $F$ be a finite field of order $q=|F|$
      [so that $q$ is a prime power and $F \simeq GF(q)$],
      let $B$ be a matrix representing $M$ over $F$,
      and let $\scrm = F^k$ for some integer $k \ge 0$.
      Then
\begin{equation}
     \Ztilde^{\rm Potts}_{F,F^k,B}(\bv)  \;=\; \Ztilde_M(q^k,\bv)
     \;.
 \label{eq.gen.FK.identity2}
\end{equation}
\end{itemize}
\end{theorem}

\noindent
We recall that an integer matrix $B$ is called
\defword{totally unimodular} if all its subdeterminants lie in $\{0,1,-1\}$.

\begin{proof}
In \reff{def.gen.Potts}, expand out the product over $e \in E$,
and let $A \subseteq E$ be the set of edges for which the term
$v_e \delta(\cdots)$ is taken.
We then need to count the number of configurations
$\bsigma = \{ \sigma_i \}_{i \in V}$ satisfying
$\sum_{i \in V} \sigma_i b_{ie} = 0$ for all $e \in A$.
In other words, we need to count the number of solutions
$\bsigma \in \scrm^V$ of the linear system $\bsigma \widehat{B} = 0$,
where $\widehat{B}$ is the submatrix of $B$ consisting of the columns from $A$.

Since both $R = \Z$ and $R=F$ are principal ideal domains,
the matrix $\widehat{B}$ has a Smith normal form
\cite{Marcus_78,Adkins_92,Brown_93,Newman_97}
\begin{equation}
   S \;=\;  P \widehat{B} Q  \;=\;  \left(\!  \begin{array}{cc}
                                                  D  & 0 \\
                                                  0  & 0
                                               \end{array}
                                    \!\right)
\end{equation}
where $P \in GL(m,R)$, $Q \in GL(n,R)$
and $D=\diag(s_1,\ldots,s_r)$;
here $m=|V|$, $n=|A|$, $r$ is the determinantal rank of $\widehat{B}$,
and $s_1,\ldots,s_r \neq 0$ are the invariant factors of $\widehat{B}$.
In case (a), since $B$ is totally unimodular,
we must have $s_1,\ldots,s_r \in \{1,-1\}$;
and by a suitable choice of $P$ (or $Q$)
we can arrange to have $s_1 = \ldots = s_r = 1$.
In case (b), because $R$ is a field,
we can again arrange to have $s_1 = \ldots = s_r = 1$.
Therefore, a vector $\bsigma \in \scrm^V$ solves $\bsigma \widehat{B} = 0$
if and only if $\bsigma P^{-1}$ is vanishing in its first $r$ components.
The number of such solutions is thus $|\scrm|^{|V|-r}$.
On the other hand, since $R$ is an integral domain,
the determinantal rank $r$ is also equal \cite[pp.~205--209]{Adkins_92}
to the column rank of $\widehat{B}$ in the quotient field of $R$
[namely, $\Q$ in case (a), and $F$ itself in case (b)];
and this rank is by definition $r_M(A)$.
This proves
\reff{eq.gen.FK.identity1}/\reff{eq.gen.FK.identity2}.\rule{.5cm}{0cm}

\vspace*{-.95cm}
\end{proof}

Theorem~\ref{thm.gen.FK}(a) implies, in particular, the known fact
\cite[Exercise~6.53(e,f)]{Brylawski_92}
that the chromatic polynomial of a regular matroid
has the consecutive-root property.
See also \cite[Theorem~III]{Crapo_69} \cite[Theorem~12.4]{Brylawski_72}
for results related to Theorem~\ref{thm.gen.FK}(a),
and \cite[Theorem~16.1]{Crapo_70} \cite[Theorem~7.6.1]{Zaslavsky_87}
for results related to Theorem~\ref{thm.gen.FK}(b).

\subsection{\mbox{\protect\boldmath $q \to 0$} limits}

The $q \to 0$ limits for the
matroid multivariate Tutte polynomial \reff{eq.defZ.matroid}
follow closely the pattern for graphs, with the obvious replacements
\begin{eqnarray*}
   \mbox{maximally connected spanning subgraph} & \longrightarrow &
   \mbox{spanning set}  \\
   \mbox{spanning forest} & \longrightarrow & \mbox{independent set}  \\
   \mbox{maximal spanning forest} & \longrightarrow & \mbox{basis}
\end{eqnarray*}

Thus, the limit $q \to 0$ with fixed $\bv$
selects out the subsets $A \subseteq E$ of maximum rank,
i.e.\ the spanning sets.
We therefore have
\begin{equation}
   \lim_{q \to 0} q^{r_M(E)} \Ztilde_M(q,\bv) \;=\;  S_M(\bv)
   \;,
\end{equation}
where $r_M(E)$ is the rank of $M$, and
\begin{equation}
   S_M(\bv) \;=\; \!\!\! \sum\limits_{\begin{scarray}
                                       A \subseteq E \\
                                       r_M(A) = r_M(E)
                                    \end{scarray}}
                       \prod_{e \in A}  v_e
\end{equation}
is the generating polynomial of \defword{spanning sets}\/ in $M$.

Similarly, the limit $q \to 0$ with fixed $\bw =\bv/q$
selects out the subsets $A \subseteq E$
having the smallest value of $|A| - r_M(A)$;
the minimum achievable value is 0,
and the sets attaining it are the independent sets.
We therefore have
\begin{equation}
   \lim_{q \to 0} \Ztilde_M(q,q\bw) \;=\;  I_M(\bw)
   \;,
\end{equation}
where
\begin{equation}
   I_M(\bw) \;=\; \!\!\! \sum\limits_{\begin{scarray}
                                       A \subseteq E \\
                                       r_M(A) = |A|
                                    \end{scarray}}
                       \prod_{e \in A}  w_e
\end{equation}
is the generating polynomial of \defword{independent sets}\/ in $M$.

Finally, we have
\begin{subeqnarray}
   \lim_{\lambda \to 0} \lambda^{-r_M(E)} S_M(\lambda\bv)  & = &  B_M(\bv) \\
   \lim_{\lambda \to \infty} \lambda^{-r_M(E)} I_M(\lambda\bw)  & = &  B_M(\bw)
\end{subeqnarray}
where
\begin{equation}
   B_M(\bv) \;=\; \!\!\!  \sum\limits_{\begin{scarray}
                                       A \subseteq E \\
                                       |A| = r_M(A) = r_M(E)
                                    \end{scarray}}
                       \prod_{e \in A}  v_e
\end{equation}
is the generating polynomial of \defword{bases}\/ in $M$.

In summary, we have the following scheme for the $q \to 0$ limits
of the multivariate Tutte polynomial for matroids:
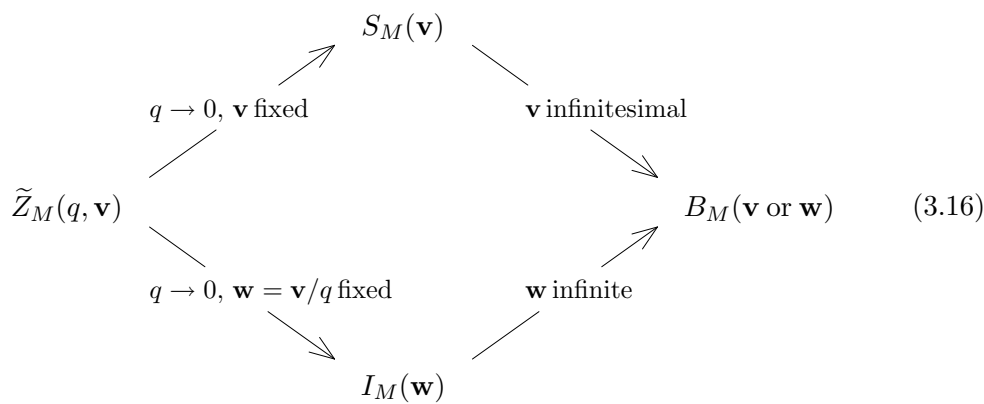
\begin{equation}
  \begin{array}{ccccc}
        &   &  S_M(\bv)   &     &    \\  
        &   \setlength{\unitlength}{1pt}
            \begin{picture}(70,50)(0,0)
            \drawline(0,0)(25.2,18)
            \put(0,23){{\small $q \to 0,\, \bv \, {\rm fixed}$}}
            \drawline(49,35)(70,50)
            \drawline(70,50)(64.34,41.76)
            \drawline(70,50)(60.37,47.32)
            \end{picture} &
        &   \setlength{\unitlength}{1pt}
            \begin{picture}(70,50)(0,0)
            \drawline(0,50)(21,35)
            \put(20,23){{\small $\bv \, {\rm infinitesimal}$}}
            \drawline(44.8,18)(70,0)
            \drawline(70,0)(60.37,2.68)
            \drawline(70,0)(64.34,8.24)
            \end{picture} &              \\
   \Ztilde_M(q,\bv)  &    &      &    &   B_M(\bv \: {\rm or} \: \bw) \\
        &   \setlength{\unitlength}{1pt}
            \begin{picture}(70,50)(0,0)
            \drawline(0,50)(21,35)
            \put(0,23){{\small $q \to 0,\, \bw = \bv/q \, {\rm fixed}$}}
            \drawline(44.8,18)(70,0)
            \drawline(70,0)(60.37,2.68)
            \drawline(70,0)(64.34,8.24)
            \end{picture} &
        &   \setlength{\unitlength}{1pt}
            \begin{picture}(70,50)(0,0)
            \drawline(0,0)(25.2,18)
            \put(20,23){{\small $\bw \, {\rm infinite}$}}
            \drawline(49,35)(70,50)
            \drawline(70,50)(64.34,41.76)
            \drawline(70,50)(60.37,47.32)
            \end{picture} & \\
        &   &  I_M(\bw)   &     &
  \end{array}
 \label{scheme.qto0.matroid}
\end{equation}

\subsection{A final remark}

Several alternative definitions of multivariate (or ``weighted'')
Tutte polynomials for graphs and/or matroids
can be found in the literature.
Most of these are essentially equivalent to $\Ztilde_M(q,\bv)$
after suitable changes of variables.
Exceptions can, however, be found in the work of Zaslavsky \cite{Zaslavsky_92}
and Bollob\'as--Riordan \cite{Bollobas_99};
their results are generalized and clarified
in an illuminating recent paper of
Ellis-Monaghan and Traldi \cite{Ellis_04}.

\section{Elementary identities}

I~now wish to prove some elementary identities
for the multivariate Tutte polynomial.
There are two alternative approaches to proving such identities:
one is to prove the identity directly for indeterminate (or complex) $q$,
using the subgraph expansion \reff{eq1.1}
or its generalization \reff{eq.defZ.matroid} to matroids;
the other is to prove the identity first for {\em positive integer}\/ $q$,
using the spin/coloring representation \reff{def.ZPotts},
and then to extend it to general $q$ by arguing that two polynomials
(or rational functions) that coincide at infinitely many points must be equal.
The latter approach is perhaps less elegant,
but it is often simpler or more intuitive.
However, only the former approach extends to matroids.

One way to guess (albeit not to prove) an identity for matroids
is to prove it first for graphs,
and then translate it from $Z_G$ to $\Ztilde_G = q^{-|V|} Z_G$;
usually the latter identity carries over verbatim to matroids,
{\em mutatis mutandis}\/.

\subsection{Disjoint unions and direct sums}

If $G$ is the disjoint union of $G_1$ and $G_2$, then trivially
\begin{equation}
   Z_G(q,\bv)  \;=\;  Z_{G_1}(q,\bv) \, Z_{G_2}(q,\bv)
   \;.
 \label{eq.components}
\end{equation}
That is, $Z_G$ ``factorizes over components''.

A slightly less trivial situation arises when $G$ consists of
subgraphs $G_1$ and $G_2$ joined at a single cut vertex $x$;
in this case
\begin{equation}
   Z_G(q,\bv)  \;=\;  {Z_{G_1}(q,\bv) \, Z_{G_2}(q,\bv)  \over q}
   \;.
 \label{eq.blocks}
\end{equation}
This is easily seen from the subgraph expansion in the form \reff{eq1.1cycles}.
It is also easily seen from the coloring representation \reff{def.ZPotts},
by first fixing the color $\sigma_x$ at the cut vertex
and then summing over it;
from this viewpoint, \reff{eq.blocks} reflects
the $S_q$ permutation symmetry of the $q$-state Potts model.\footnote{
   More precisely, it reflects the symmetry of the spin model
   under a global transformation $\sigma_y \mapsto g\sigma_y$
   (simultaneously for all $y \in V$)
   that acts {\em transitively}\/ on each single-spin space.
}
We summarize \reff{eq.blocks} by saying that
$Z_G$ ``factorizes over blocks'' modulo a factor $q$.

The identities \reff{eq.components} and \reff{eq.blocks}
can be written in a unified form, by using $\Ztilde_G = q^{-|V|} Z_G$:
in both cases we have
\begin{equation}
   \Ztilde_G(q,\bv)  \;=\;  \Ztilde_{G_1}(q,\bv) \, \Ztilde_{G_2}(q,\bv)
   \;.
\end{equation}
This, in turn, is a special case of the following obvious fact:
if a matroid $M$ is the direct sum of matroids $M_1$ and $M_2$, then
\begin{equation}
   \Ztilde_M(q,\bv)  \;=\;  \Ztilde_{M_1}(q,\bv) \, \Ztilde_{M_2}(q,\bv)
   \;.
\end{equation}

\subsection{Duality}

Suppose first that $G=(V,E)$ is a connected {\em planar}\/ graph.
Consider any plane embedding of $G$,
and let $G^* = (V^*,E^*)$ be the corresponding dual graph.
There is a natural bijection between $E$ and $E^*$
(namely, an edge $e \in E$ is identified with the unique edge
$e^* \in E^*$ that it crosses),
so we shall henceforth identify $E^*$ with $E$.
Of course, the vertex set $V^*$ can be identified with the faces
in the given embedding of $G$, so by Euler's relation we have
\begin{equation}
   |V| - |E| + |V^*|  \;=\;  2  \;.
 \label{eq.euler}
\end{equation}

Consider now any subset $A \subseteq E$,
and draw in $G^*$ the {\em complementary}\/ set of edges ($E \setminus A$).
Simple topological arguments then yield the relations
\begin{subeqnarray}
   k_{G}(A)                & = &  c_{G^*}(E \setminus A)  \,+\,  1  \\[2mm]
   k_{G^*}(E \setminus A)  & = &  c_{G}(A)  \,+\,  1
  \label{eq.dualranks}
\end{subeqnarray}
where $k=$ components and $c=$ cyclomatic number.
[Note that (\ref{eq.dualranks}a) and (\ref{eq.dualranks}b)
 are equivalent, as a consequence of \reff{eq.euler}.]
Substituting \reff{eq.dualranks} into \reff{eq2.1}--\reff{eq1.1cycles},
we deduce the \defword{duality relation}\/ \cite{Wu_82}
\begin{equation}
   Z_{G^*}(q, \bv)  \;=\;
   q^{1-|V|} \left( \prod\limits_{e \in E} v_e \right)  Z_G(q, q/\bv)
   \;.
 \label{eq.duality}
\end{equation}
In brief, duality takes $v_e \mapsto q/v_e$
(and inserts some prefactors).
Note that two applications of \reff{eq.duality}
lead us back to where we started, thanks to \reff{eq.euler}.
In terms of $\Ztilde_G = q^{-|V|} Z_G$, we have
\begin{subeqnarray}
   \Ztilde_{G^*}(q, \bv)  & = &
   q^{1-|V^*|} \left( \prod\limits_{e \in E} v_e \right)  \Ztilde_G(q, q/\bv)
       \\[2mm]
   & = &
   q^{|V|-1} \left( \prod\limits_{e \in E} {v_e \over q} \right)
                                                          \Ztilde_G(q, q/\bv)
   \;.
 \label{eq.duality.Ztilde}
\end{subeqnarray}

In the $q \to 0$ limit, we obtain the following
special cases of \reff{eq.duality}:
\begin{eqnarray}
   C_{G^*}(\bv)   & = &   \left( \prod\limits_{e \in E} v_e \right)
                                  F_G(1/\bv)
     \label{eq.duality1} \\[2mm]
   F_{G^*}(\bw)   & = &   \left( \prod\limits_{e \in E} w_e \right)
                                  C_G(1/\bw)
     \label{eq.duality2} \\[2mm]
   T_{G^*}(\bv)   & = &   \left( \prod\limits_{e \in E} v_e \right)
                                  T_G(1/\bv)
     \label{eq.duality3}
\end{eqnarray}
We can also relate the multivariate Tutte polynomial on $G^*$
to the multivariate flow polynomial on $G$:
by \reff{thm.potts.flowpoly.eq1} we have
\begin{equation}
   Z_{G^*}(q, \bv)  \;=\;  q F_G(q, \bv)   \;.
\end{equation}
In particular, the chromatic polynomial on $G^*$
is essentially identical to the flow polynomial on $G$:
\begin{equation}
   P_{G^*}(q)  \;=\;  q F_G(q)  \;.
\end{equation}
The chromatic polynomial and the flow polynomial
are thus ``dual'' objects.

Among graphs, only {\em planar}\/ graphs have duals with good properties
\cite[section~4.6]{Diestel_97} \cite[section~5.2]{Oxley_92};
one major advantage of considering matroids
is that {\em every}\/ matroid has a dual.
The duality formula for the matroid multivariate Tutte polynomial is,
not surprisingly, identical in form to \reff{eq.duality.Ztilde}:
\begin{subeqnarray}
   \Ztilde_{M^*}(q, \bv)  & = &
   q^{-r_{M^*}(E)} \left( \prod\limits_{e \in E} v_e \right)
                                                          \Ztilde_M(q, q/\bv)
       \\[2mm]
   & = &
   q^{r_M(E)} \left( \prod\limits_{e \in E} {v_e \over q} \right)
                                                          \Ztilde_M(q, q/\bv)
   \;.
 \label{eq.duality.matroid}
\end{subeqnarray}
[Here $r_M(E)$ is the rank of $M$ and $r_{M^*}(E)$ is the rank of $M^*$;
 their sum is $|E|$.]
Indeed, \reff{eq.duality.matroid} is an easy consequence of the
definition \reff{eq.defZ.matroid} together with the
formula for the rank function of a dual:
\begin{equation}
   r_{M^*}(A)  \;=\;  |A| \,+\, r_M(E \setminus A) \,-\, r_M(E)
   \;.
\end{equation}

\subsection{Deletion-contraction identity}

If $e \in E$, let $G \setminus e$
denote the graph obtained from $G$ by deleting the edge $e$,
and let $G/e$ denote the graph obtained from $G \setminus e$
by contracting the two endpoints of $e$ into a single vertex
(please note that we retain in $G/e$ any loops or multiple edges
 that may be formed as a result of the contraction).
Then, for any $e \in E$, we have the identity
\begin{equation}
   Z_G(\qvbf)   \;=\;
      Z_{G \setminus e}(q, \bv_{\neq e})   \,+\,
      v_e Z_{G/e}(q, \bv_{\neq e})
 \label{eq.delcon}
\end{equation}
where $\bv_{\neq e} = \{ v_f \} _{f \in E \setminus e}$.
This is easily seen either from the coloring representation \reff{def.ZPotts}
or the subgraph expansion \reff{eq1.1}.
Please note that the deletion-contraction identity \reff{eq.delcon}
takes the same form regardless of whether $e$ is a normal edge,
a loop, or a bridge
(in contrast to the situation for the usual Tutte polynomial $T_G$).
Of course, if $e$ is a loop, then $G/e = G \setminus e$,
so we can also write $Z_G = (1+v_e) Z_{G \setminus e} = (1+v_e) Z_{G/e}$.
Similarly, if $e$ is a bridge, then $G \setminus e$ is the disjoint union
of two subgraphs $G_1$ and $G_2$ while $G/e$ is obtained by joining
$G_1$ and $G_2$ at a cut vertex,
so that $Z_{G/e} = Z_{G \setminus e}/q$ and hence
$Z_G = (1+v_e/q) Z_{G \setminus e} = (q+v_e) Z_{G/e}$.


In terms of $\Ztilde_G = q^{-|V|} Z_G$,
the deletion-contraction identity takes the form
\begin{subeqnarray}
   \Ztilde_G  & = &  \Ztilde_{G \setminus e} \,+\, {v_e \over q} \Ztilde_{G/e}
                        \qquad\hbox{if $e$ is not a loop}   \\[3mm] 
   \Ztilde_G  & = &  \Ztilde_{G \setminus e} \,+\, v_e \Ztilde_{G/e}
                                                       \nonumber \\
              & = &  (1+v_e) \Ztilde_{G \setminus e}   \nonumber \\
              & = &  (1+v_e) \Ztilde_{G/e}
                        \qquad\quad\;\hbox{if $e$ is a loop}
 \label{eq.delcon.Ztilde}
\end{subeqnarray}
as easily follows from \reff{eq.delcon}
together with the counting of vertices in $G \setminus e$ and $G/e$.

Not surprisingly, the deletion-contraction formula for matroids
is identical in form to \reff{eq.delcon.Ztilde}:
\begin{subeqnarray}
   \Ztilde_M  & = &  \Ztilde_{M \setminus e} \,+\, {v_e \over q} \Ztilde_{M/e}
                        \qquad\hbox{if $e$ is not a loop}   \\[3mm] 
   \Ztilde_M  & = &  \Ztilde_{M \setminus e} \,+\, v_e \Ztilde_{M/e}
                                                       \nonumber \\
              & = &  (1+v_e) \Ztilde_{M \setminus e}   \nonumber \\
              & = &  (1+v_e) \Ztilde_{M/e}
                        \qquad\quad\;\hbox{if $e$ is a loop}
 \label{eq.delcon.matroid}
\end{subeqnarray}
This easily follows from the formulae for the rank function
of a deletion or contraction:  if $A \subseteq E \setminus e$, then
\begin{subeqnarray}
   r_{M \setminus e}(A)   & = &   r_M(A)   \\[3mm]
   r_{M/e}(A)             & = &
        \begin{cases}
            r_M(A \cup e) - 1  & \text{if $e$ is not a loop} \\
            r_M(A \cup e)      & \text{if $e$ is a loop}
        \end{cases}
  \slabel{eq.rank.contraction}
\end{subeqnarray}

{\bf Remark.}
Many treatments in the literature {\em define}\/ the Tutte polynomial
by the deletion-contraction identity \reff{eq.delcon}
together with the initial condition $Z_G(q, {\bf v}) = q^{|V|}$
for an edgeless graph
(or the analogous thing for the standard Tutte polynomial $T_G$).
But this approach has the disadvantage that one must prove
that this $Z_G$ is well-defined,
i.e.\ that the result does not depend on the order in which
one applies \reff{eq.delcon} to the various edges.
A much cleaner approach, it seems to me, is to define $Z_G(q, {\bf v})$
by the explicit formula \reff{eq1.1},
and then deduce the deletion-contraction identity as an immediate property.

\subsection{Parallel-reduction identity}  \label{subsec.parallel}

If $G$ contains edges $e_1,e_2$
connecting the same pair of vertices $x,y$, they can be replaced,
without changing the value of $Z$,
by a single edge $e=xy$ with weight
\begin{equation}
   v_e  \;=\;  (1+v_{e_1})(1+v_{e_2}) \,-\, 1
        \;=\;  v_{e_1} + v_{e_2} + v_{e_1} v_{e_2}
   \;.
 \label{eq.parallel1}
\end{equation}
This is easily seen either from the coloring representation \reff{def.ZPotts}
or the subgraph expansion \reff{eq1.1}.
More formally, we can identify the new edge $e$
with (for instance) the old edge $e_1$ after deletion of $e_2$,
and thus write
\begin{equation}
   Z_G(q, \bv_{\neq e_1,e_2}, v_{e_1}, v_{e_2})  \;=\;
   Z_{G \setminus e_2}(q, \bv_{\neq e_1,e_2},
                          v_{e_1} + v_{e_2} + v_{e_1} v_{e_2})
   \;.
 \label{eq.parallel2}
\end{equation}
The parallel-reduction rule
$(v_1,v_2) \mapsto v_{\rm eff}$ with $1 + v_{\rm eff} = (1+v_1)(1+v_2)$
can be remembered by the mnemonic ``$1+v$ multiplies''.
We write $v_1 \parallel v_2 \equiv (1+v_1)(1+v_2) - 1$.


A virtually identical formula holds for matroids:
if $e_1$ and $e_2$ are parallel elements in a matroid $M$
(i.e., form a two-element circuit), then
\begin{equation}
   \Ztilde_M(q, \bv_{\neq e_1,e_2}, v_{e_1}, v_{e_2})  \;=\;
   \Ztilde_{M \setminus e_2}(q, \bv_{\neq e_1,e_2},
                          v_{e_1} + v_{e_2} + v_{e_1} v_{e_2})
   \;.
 \label{eq.parallel3}
\end{equation}
The formula \reff{eq.parallel3} also holds trivially if
$e_1$ and $e_2$ are both loops.

\subsection{Series-reduction identity}  \label{subsec.series}

We say that edges $e_1, e_2 \in E$
are {\em in series (in the narrow sense)}\/
if there exist vertices $x,y,z \in V$ with $x \neq y$ and $y \neq z$
such that $e_1$ connects $x$ and $y$,
$e_2$ connects $y$ and $z$, and $y$ has degree 2 in $G$.
In this case the pair of edges $e_1,e_2$ can be replaced,
without changing the value of $Z$,
by a single edge $e=xz$ with weight
\begin{equation}
   v_e  \;=\;  {v_{e_1} v_{e_2}  \over  q + v_{e_1} + v_{e_2}}
 \label{eq.series1}
\end{equation}
provided that we then multiply $Z$ by the prefactor $q + v_{e_1} + v_{e_2}$.
More formally, we can identify the new edge $e$
with (for instance) the old edge $e_1$ after contraction of $e_2$,
and thus write
\begin{equation}
   Z_G(q, \bv_{\neq e_1,e_2}, v_{e_1}, v_{e_2})  \;=\;
   (q + v_{e_1} + v_{e_2}) \,
   Z_{G / e_2}(q, \bv_{\neq e_1,e_2},
                          v_{e_1} v_{e_2} / (q + v_{e_1} + v_{e_2}))
   \;.
 \label{eq.series2}
\end{equation}
This identity can be derived from the coloring representation \reff{def.ZPotts}
by noting that
\begin{subeqnarray}
 & &
   \sum_{\sigma_y = 1}^q
   \bigl[ 1 + v_{e_1} \delta(\sigma_x,\sigma_y) \bigr]
   \,
   \bigl[ 1 + v_{e_2} \delta(\sigma_y,\sigma_z) \bigr]
 \nonumber \\
 & & \qquad =\;
   q + v_{e_1} + v_{e_2} + v_{e_1} v_{e_2} \delta(\sigma_x,\sigma_z)  \\
 & & \qquad =\;
   (q + v_{e_1} + v_{e_2})
   \left[1 \,+\, {v_{e_1} v_{e_2}  \over  q + v_{e_1} + v_{e_2}}
                 \delta(\sigma_x,\sigma_z)
   \right]
   \;.
\end{subeqnarray}
Alternatively, it can be derived from the subgraph expansion \reff{eq1.1}
by considering the four possibilities for the edges $e_1$ and $e_2$
to be occupied or empty and analyzing the number of connected components
thereby created.
The series-reduction rule
$(v_1,v_2) \mapsto v_{\rm eff} \equiv v_1 v_2/(q+v_1+v_2)$
can be remembered by the mnemonic ``$1+q/v$ multiplies'':
namely,
\begin{equation}
   1 + {q \over v_{\rm eff}}  \;=\;
   \left( 1 + {q \over v_1} \right)
   \left( 1 + {q \over v_2} \right)
   \;.
\end{equation}
We write $v_1 \seriesq v_2 \equiv v_1 v_2/(q+v_1+v_2)$.

Consider now the more general situation in which
$\{e_1,e_2\}$ is a two-edge cut of $G$
(not necessarily the cut associated with a degree-2 vertex $y$);
we then say that $e_1,e_2$ are {\em in series (in the wide sense)}\/.
It turns out that the identity \reff{eq.series2} still holds.
To see this, let us prove the generalization of this identity to matroids.
Let $e_1$ and $e_2$ be series elements in a matroid $M$,
i.e., suppose that $\{e_1,e_2\}$ is a cocircuit.
Then, for any $A \subseteq E \setminus \{e_1,e_2\}$, we have
\begin{equation}
   r_M(A \cup e_1)  \;=\; r_M(A \cup e_2)  \;=\;  r_M(A) + 1
 \label{eq.rank.cocircuit}
\end{equation}
(since the complement of a cocircuit is a hyperplane).
A short calculation using \reff{eq.rank.contraction} with $e=e_2$ then yields
\begin{equation}
   \Ztilde_M(q, \bv_{\neq e_1,e_2}, v_{e_1}, v_{e_2})  \;=\;
   {q + v_{e_1} + v_{e_2} \over  q}
   \,
   \Ztilde_{M / e_2}(q, \bv_{\neq e_1,e_2},
                          v_{e_1} v_{e_2} / (q + v_{e_1} + v_{e_2}))
   \;.
 \label{eq.series3}
\end{equation}
%
%
The formula \reff{eq.series3} also holds trivially if
$e_1$ and $e_2$ are both coloops.

Please note that duality $v \mapsto q/v$ interchanges
the parallel-reduction rule (``$1+v$ multiplies'')
with the series-reduction rule (``$1+q/v$ multiplies'').
This is no accident,
since we now see that parallel-reduction and series-reduction
(in the wide sense)
are indeed duals of each other:
$\{e_1,e_2\}$ is a circuit (resp.\ cocircuit) in $M$
if and only if it is a cocircuit (resp.\ circuit) in $M^*$.

{\bf Remark.}
Using the series-reduction formula \reff{eq.series2}
together with the multivariate approach,
one can give a simple proof of (a generalization of)
the Brown--Hickman \cite{Brown-Hickman_99b} theorem on chromatic roots
of large subdivisions \cite[Appendix A]{Sokal_chromatic_roots}.
Likewise, using the parallel- and series-reduction formulae
\reff{eq.parallel2}/\reff{eq.series2},
one can give a simple proof of (a slight generalization of)
Thomassen's \cite{Thomassen_97} construction of
2-degenerate graphs with arbitrarily large real chromatic roots
\cite[Appendix B]{Sokal_chromatic_roots}.

\subsection{Reduction formulae for 2-rooted subgraphs}
   \label{subsec.2-rooted}

Let $G=(V,E)$ be a finite graph, and let $x,y$ be distinct vertices of $G$.
We define $G/xy$ to be the graph in which $x$ and $y$ are contracted
to a single vertex.  (N.B.:  If $G$ contains one or more edges $xy$,
then these edges are {\em not}\/ deleted, but become loops in $G/xy$.)
There is a canonical one-to-one correspondence between the edges of $G$
and the edges of $G/xy$;  for simplicity (though by slight abuse of notation)
we denote an edge of $G$ and the corresponding edge of $G/xy$
by the same letter.  In particular, we can apply a given set of edge weights
$\{ v_e \} _{e \in E}$ to both $G$ and $G/xy$.

Let us now define
\begin{eqnarray}
   \zgxconny(\qvbf)   & = & \!\!\!
   \sum_{\begin{scarray}
            A \subseteq E \\
            A \, {\rm connects} \, x \, {\rm to} \, y
         \end{scarray}
        }
   \!\!\!\!\!  q^{k(A)} \;  \prod_{e \in A}  v_e
   \\[4mm]
   \zgxnoconny(\qvbf)   & = &
   \!\!\!
   \sum_{\begin{scarray}
            A \subseteq E \\
            A \, {\rm does\,not\,connect} \, x \, {\rm to} \, y
         \end{scarray}
        }
   \!\!\!\!\!  q^{k(A)} \;  \prod_{e \in A}  v_e
\end{eqnarray}
{}From \reff{eq1.1} we have trivially
\begin{equation}
   Z_G(\qvbf)  \;=\;  \zgxconny(\qvbf) \,+\, \zgxnoconny(\qvbf)
 \label{eq2.G}
\end{equation}
and almost as trivially
\begin{equation}
   Z_{G/xy}(\qvbf)  \;=\;  \zgxconny(\qvbf) \,+\, q^{-1} \zgxnoconny(\qvbf)
   \;.
 \label{eq2.Gxy}
\end{equation}

%

Let now $q$ be an integer $\ge 1$,
and define the restricted Potts-model partition function
\begin{equation}
   Z^{\rm Potts}_{G,x,y}(q, \bv; \sigma_x, \sigma_y)   \;=\;
   \sum_{ \sigma\colon\, V \setminus \{x,y\} \to \{1,\ldots,q\}  }
   \,  \prod_{e \in E}  \,
      \biggl[ 1 + v_e \delta(\sigma_{x_1(e)}, \sigma_{x_2(e)}) \biggr]
 \label{eq.ZPotts.rest}
\end{equation}
where $\sigma_x, \sigma_y \in \{1,\ldots,q\}$.
We then have the following refinement of the
Fortuin--Kasteleyn identity \reff{eq.FK.identity}:

\begin{prop}
  \label{prop2.1}
\begin{equation}
   Z^{\rm Potts}_{G,x,y}(q, \bv; \sigma_x, \sigma_y)   \;=\;
   A_{G,x,y}(q, \bv)  \,+\,
      B_{G,x,y}(q, \bv) \, \delta(\sigma_x,\sigma_y)
 \label{eq_prop2.1_1}
\end{equation}
where
\begin{subeqnarray}
   A_{G,x,y}(q, \bv)   & = &  q^{-2} \zgxnoconny(\qvbf)  \\[2mm]
   B_{G,x,y}(q, \bv)   & = &  q^{-1} \zgxconny(\qvbf)
 \label{eq_prop2.1_2}
\end{subeqnarray}
are polynomials in $q$ and $\{v_e\}$, whose degrees in $q$ are
\begin{subeqnarray}
   \deg A   & = &  |V| \,-\, 2  \\
   \deg B   &\le&  |V| \,-\, 1 \,-\, {\rm dist}(x,y)
\end{subeqnarray}
where ${\rm dist}(x,y)$ is the length of the shortest path
from $x$ to $y$ using edges having $v_e \neq 0$
(if no such path exists, then $B = 0$).
\end{prop}

\begin{firstproof}
In \reff{eq.ZPotts.rest}, expand out the product over $e \in E$,
and let $A \subseteq E$ be the set of edges for which the term
$v_e \delta(\sigma_{x_1(e)}, \sigma_{x_2(e)})$ is taken.
Now perform the sum over configurations
$\{ \sigma_z \} _{z \in V \setminus \{x,y\}}$:
in each connected component of the subgraph $(V,A)$
the spin value $\sigma_z$ must be constant.
In particular, in each component containing $x$ and/or $y$,
the spins must all equal the specified value $\sigma_x$ and/or $\sigma_y$;
in all other components, the spin value is free.
Therefore,
\begin{equation}
   Z^{\rm Potts}_{G,x,y}(q, \bv; \sigma_x, \sigma_y)   \;=\;
   \begin{cases}
         q^{-2} \zgxnoconny(\qvbf)   & \text{if $\sigma_x \neq \sigma_y$} \\
         q^{-2} \zgxnoconny(\qvbf) + q^{-1} \zgxconny(\qvbf)
                                     & \text{if $\sigma_x = \sigma_y$}
   \end{cases}
\end{equation}
The claims about the degrees are easily verified.
\end{firstproof}

\begin{secondproof}
For each positive integer $q$,
the $S_q$ permutation symmetry of the Potts model implies that
\begin{equation}
   Z^{\rm Potts}_{G,x,y}(q, \bv; \sigma_x, \sigma_y)   \;=\;
      A \,+\, B \delta(\sigma_x,\sigma_y)
\end{equation}
for some numbers $A$ and $B$ depending on $G$, $x$, $y$, $q$ and $\bv$.
But then, summing over $\sigma_x, \sigma_y$ without and with the constraint
$\sigma_x = \sigma_y$, we get
\begin{subeqnarray}
   Z_G                & = &  q^2 A \,+\, q B   \\
   Z_{G/xy}   & = &  q A \,+\, q B
\end{subeqnarray}
Hence
\begin{subeqnarray}
   A   & = &   {Z_G - Z_{G/xy}   \over   q(q-1)}
       \;=\;   q^{-2} \zgxnoconny
   \\[2mm]
   B   & = &   {q Z_{G/xy} - Z_G   \over   q(q-1)}
       \;=\;   q^{-1} \zgxconny
 \label{eq2.AB}
\end{subeqnarray}
by virtue of \reff{eq2.G} and \reff{eq2.Gxy}.
\end{secondproof}

Let us now consider inserting the 2-rooted graph $(G,x,y)$
in place of an edge $e_*$ in some other graph $H$.
{}From \reff{eq_prop2.1_1}/\reff{eq_prop2.1_2}
we see that $(G,x,y)$ then acts as a single edge with effective weight
\begin{equation}
   v_{\rm eff}  \;\equiv\;  v_{{\rm eff},G,x,y}(\qvbf)
    \;=\;  {B_{G,x,y}(\qvbf) \over A_{G,x,y}(\qvbf)}
    \;=\;  {q \zgxconny(\qvbf) \over \zgxnoconny(\qvbf)}
   \;,
\end{equation}
which is a rational function of $q$ and $\{v_e\}$;
in addition, the partition function is multiplied by
an overall prefactor $A_{G,x,y}(\qvbf)$.
This follows from Proposition~\ref{prop2.1}
whenever $q$ is an integer $\ge 1$;
and the corresponding identity then holds for all $q$,
because both sides are rational functions of $q$
that agree at infinitely many points.
It is also worth noting that the ``transmissivity''
$t_{{\rm eff}} \equiv v_{{\rm eff}} / (q + v_{{\rm eff}})$
is given by the simple formula
\begin{equation}
   t_{{\rm eff}}
    \;=\;  {\zgxconny(\qvbf) \over Z_G(\qvbf)}
   \;.
\end{equation}

The most general version of this construction appears to be the following
\cite{Woodall_00}:
Let $H=(V,E)$ be a finite undirected graph,
and let $\vec{H}$ be a directed graph obtained by assigning an orientation
to each edge of $H$.
For each edge $e \in E$, let $G_e = (V_e, E_e, x_e, y_e)$
be a 2-rooted finite undirected graph
(so that $x_e, y_e \in V_e$ with $x_e \neq y_e$)
equipped with edge weights $\{ v_{\widetilde{e}} \} _{\widetilde{e} \in E_e}$.
We denote by ${\bf G}$ the family $\{ G_e \} _{e \in E}$,
and we denote by $\vec{H}^{\bf G}$ the undirected graph obtained from $H$
by replacing each edge $e \in E$ with a copy of the corresponding graph $G_e$,
attaching $x_e$ to the tail of $e$ and $y_e$ to the head.
Its edge set is thus
${\bf E} = \biguplus\limits_{e \in E} E_e$ (disjoint union).
We then have:

\begin{prop}
   \label{prop2.2}
Let $H=(V,E)$ and $\{ G_e \} _{e \in E}$ be as above.
Suppose that
\begin{eqnarray}
& &
   Z_{G_e,x_e,y_e}(q, \{v_{\widetilde{e}} \} _{\widetilde{e} \in E_e};
                   \sigma_{x_e}, \sigma_{y_e})
   \nonumber \\
& & \qquad=\;
   A_{G_e,x_e,y_e}(q, \{v_{\widetilde{e}} \} _{\widetilde{e} \in E_e})
   \,+\,
   B_{G_e,x_e,y_e}(q, \{v_{\widetilde{e}} \} _{\widetilde{e} \in E_e}) \,
      \delta(\sigma_{x_e},\sigma_{y_e})
   \;,
   \qquad
\end{eqnarray}
and define
\begin{equation}
   v_{{\rm eff},e}  \;\equiv\;
    {B_{G_e,x_e,y_e}(q, \{v_{\widetilde{e}} \} _{\widetilde{e} \in E_e})
     \over
     A_{G_e,x_e,y_e}(q, \{v_{\widetilde{e}} \} _{\widetilde{e} \in E_e})
    }
   \;.
\end{equation}
Then
\begin{equation}
   Z_{\vec{H}^{\bf G}}(q, \{v_{\widetilde{e}} \} _{\widetilde{e} \in {\bf E}})
   \;=\;
   \left( \prod\limits_{e \in E}
          A_{G_e,x_e,y_e}(q, \{v_{\widetilde{e}} \} _{\widetilde{e} \in E_e})
   \right)
   \,\times\,
   Z_H(q, \{v_{{\rm eff},e} \} _{e \in E} )
   \;.
\end{equation}
In particular, $Z_{\vec{H}^{\bf G}}$ does not depend on the orientations
of the edges of $H$.
\end{prop}

\subsection{Analogy with electrical circuit theory}

The identities discussed in
Sections~\ref{subsec.parallel}--\ref{subsec.2-rooted}
are strongly reminiscent of analogous identities
in electrical circuit theory:
there are elementary formulae for the reduction of
linear circuit elements placed in series or parallel;
and more generally, any 2-terminal subnetwork consisting of
linear passive circuit elements is equivalent
(at a fixed frequency $\omega$) to some single ``effective admittance''.
But the relation of the multivariate Tutte polynomial
to electrical circuit theory goes beyond mere analogy.
For, as was discovered by Kirchhoff \cite{Kirchhoff_1847} in 1847
and will be reviewed in Section~\ref{sec.spanning},
linear electrical circuits are intimately related to spanning trees;
and as was seen in \reff{scheme.qto0},
the spanning-tree polynomial $T_G(\bv)$ arises from the double limit of the
multivariate Tutte polynomial in which $q \to 0$ and
$\bv$ is infinitesimal.

Look back, for instance, at the Tutte--Potts parallel law
$v_1 \parallel v_2 \equiv v_1 + v_2 + v_1 v_2$;
it is a nonlinear generalization of the familiar law
$v_1 \parallel v_2 \equiv v_1 + v_2$
for combining electrical conductances in parallel,
and reduces to it when $\bv$ is infinitesimal.
Likewise, the Tutte--Potts series law
$v_1 \seriesq v_2 \equiv v_1 v_2/(q+v_1+v_2)$
is a nonlinear generalization of the familiar law
$v_1 \seriesnoq v_2 \equiv v_1 v_2/(v_1+v_2)$
for combining electrical conductances in series,
and reduces to it when $q \to 0$.

It thus makes sense to ask which results of electrical circuit theory
can be generalized from the spanning-tree polynomial $T_G(\bv)$
to the connected-spanning-subgraph polynomial $C_G(\bv)$,
the spanning-forest polynomial $F_G(\bw)$,
or even the full multivariate Tutte polynomial $Z_G(q,\bv)$.

\section{Complex zeros of \mbox{\protect\boldmath $Z_G$}: Why should we care?}
   \label{sec.complex}

Since $Z_G(q, \bv)$ is a polynomial in $q$ and $\{ v_e \}$,
it is mathematically natural to inquire about the zeros
of $Z_G(q, \bv)$ as a function of $q$ or $\{ v_e \}$ or both.
Here we can treat $q$ and $\{ v_e \}$ as either real or complex variables.

The real zeros of the chromatic polynomial have been extensively studied,
and to a lesser extent this is true also for the flow polynomial;
see \cite{Jackson_03} for an excellent recent survey.
Many of these results for the chromatic polynomial
extend to the multivariate Tutte polynomial $Z_G(q,\bv)$
in the antiferromagnetic regime (i.e., $-1 \le v_e \le 0$ for all $e$)
\cite{Jackson-Sokal_zerofree}.\footnote{
   Indeed, some of the proofs become {\em simpler}\/
   when the multivariate generalizations are considered:
   see \cite{Jackson-Sokal_zerofree} for details.
   In particular, the multivariate approach clarifies the meaning
   of the mysterious number $q=32/27$ \cite{Jackson_93,Edwards_98}.
}

The complex zeros of the partition function play
a very important role in
statistical mechanics,
as was shown by Yang and Lee \cite{Yang-Lee_52} in 1952.
This arises as follows:
Statistical physicists are interested in \defword{phase transitions}\/,
namely in points where one or more physical quantities
(e.g.\ the energy or the magnetization)
depend nonanalytically (in many cases even discontinuously)
on one or more control parameters
(e.g.\ the temperature or the magnetic field).
Now, such nonanalyticity is manifestly impossible in
\reff{eq1.1}/\reff{def.ZPotts} ---
or in any other reasonable statistical-mechanical model, for that matter ---
for any finite graph $G$.
Rather, phase transitions arise only in the
\defword{infinite-volume limit}\/.
That is, we consider some countably infinite graph
$G_\infty = (V_\infty, E_\infty)$
--- usually a regular lattice, such as $\Z^d$ with nearest-neighbor edges ---
and an increasing sequence of finite subgraphs $G_n = (V_n, E_n)$.
It can then be shown (under modest hypotheses on the $G_n$)
that the \defword{(limiting) free energy per unit volume}\/
\begin{equation}
   f_{G_\infty}(q,v)   \;=\;
   \lim_{n \to \infty}   |V_n|^{-1}  \log Z_{G_n}(q,v)
 \label{limiting_free_energy}
\end{equation}
exists for all nondegenerate physical values
of the parameters\footnote{
   Here ``physical'' means that the weights are nonnegative,
   so that the model has a probabilistic interpretation;
   and ``nondegenerate'' means that we exclude the limiting cases
   $v=-1$ in (a) and $q=0$ in (b), which cause difficulties
   due to the existence of configurations having zero weight.
},
namely either
\begin{quote}
\begin{itemize}
  \item[(a)]     $q$ integer $\ge 1$ and $-1 < v < \infty$
     \quad  (using \reff{def.ZPotts},
        see \cite[Section I.2]{Israel_79})
  \item[or (b)]  $q$ real $> 0$ and $0 \le v < \infty$
     \quad  (using \reff{eq1.1},
        see \cite[Theorem 4.1]{Grimmett_95} and
        \cite{Grimmett_78,Seppalainen_98})
\end{itemize}
\end{quote}
This limit $f_{G_\infty}(q,v)$ is in general a continuous function of $v$
[and of $q$ in case (b)];
but it can fail to be a real-analytic function of $v$ and/or $q$,
because complex singularities of $\log Z_{G_n}(q,v)$
--- namely, complex zeros of $Z_{G_n}(q,v)$ ---
can approach the real axis in the limit $n \to\infty$.
Therefore, the possible points of physical phase transitions
are precisely the real limit points of such complex zeros.
As a result, theorems that constrain the possible location of
complex zeros of the partition function are of great interest.\footnote{
   I don't want to give the impression that the partition function
   is the {\em only}\/ quantity of interest to statistical physicists.
   Indeed, what is really of physical significance is not the
   partition function, but rather the
   {\em Boltzmann--Gibbs probability measure}\/,
   which gives the probabilities for different configurations of the system:
   in the Potts (resp.\ random-cluster) model this measure is given by the
   summand of \reff{def.ZPotts} [resp.\ \reff{eq1.1bis}]
   divided by the partition function $Z_G(q,\bv)$.
   From this point of view, the partition function is merely a
   normalization factor, of no intrinsic physical interest whatsoever!
   Nevertheless, {\em derivatives}\/ of the partition function
   with respect to the parameters occurring in it yield
   correlation functions (i.e., expectations of particular random variables
   with respect to the Boltzmann--Gibbs measure),
   which {\em are}\/ of direct physical interest.
   Furthermore, the complex zeros of the partition function
   are of interest for the reason just explained;
   the Yang--Lee picture thus provides {\em one}\/ method for investigating
   phase transitions, to be used in conjunction with other approaches.
}

In particular, if it can be proven that the partition function $Z_{G_n}$
is nonvanishing in some complex domain $D$
{\em that does not vary with the size of the graph $G_n$}\/,
then (under mild hypotheses) the infinite-volume free energy $f_{G_\infty}$
will be analytic in $D$ ---
that is, in physical terms, there will be no phase transitions in $D$.
Study of the complex zeros of the partition function
is thus one way of proving the absence of phase transitions
in specified regions of parameter space.

It is useful to distinguish between ``soft'' and ``hard'' theorems
of this type.
The ``soft'' theorems assert the nonvanishing of the partition function $Z_G$
for some specified statistical-mechanical model
in some ``natural'' large domain of multidimensional complex space
(e.g.\ the unit polydisc or the product of right half-planes)
for arbitrary graphs $G$.
The ``hard'' theorems, by contrast, assert the nonvanishing of $Z_G$
in a domain that depends quantitatively on some parameter
associated to $G$, such as its maximum degree.
(``Hard'' theorems are of physical interest only when the parameter
 in question is ``local'', so that the domain obtained is uniform in
 the $G_n$.  E.g.~maximum degree is OK, but total number of vertices
 is not.)

The ``soft'' theorems are particularly interesting
from a mathematical point of view:
they assert that some combinatorially natural generating polynomial
is nonvanishing in some large domain of $\C^n$.
I~am aware of only three classes of such theorems
(but I~urge readers to try to discover new ones!):
\begin{itemize}
   \item[1)] The Lee--Yang theorem for the ferromagnetic Ising model
       at complex magnetic field \cite{Lee-Yang_52}
       and its generalizations.\footnote{
           A partial bibliography (through 1980) of generalizations
           of the Lee--Yang theorem can be found in \cite{Lieb-Sokal_81}.
}
   \item[2)]  The Heilmann--Lieb theorem for the matching polynomial
       \cite{Heilmann_72}.\footnote{
           This theorem is most often quoted in the form of its
           univariate corollary,
           namely that the roots of the univariate matching polynomial
           are all (negative) real.
           But it is the multivariate theorem that is truly fundamental.
           Two simple proofs of the multivariate Heilmann--Lieb theorem
           can be found in
           \cite[Theorem~10.1 and Remark 4 following it]{Choe_HPP}.
}
   \item[3)]  The half-plane theorem for the spanning-tree polynomial
       and its generalizations \cite{Choe_HPP}.
\end{itemize}
The Lee--Yang and Heilmann--Lieb theorems can, in fact,
be given a unified combinatorial formulation:
Let us equip the graph $G = (V,E)$ with
edge weights $\boldsymbol{\lambda} = \{\lambda_e\}_{e \in E}$
and vertex weights ${\bf t} = \{t_i\}_{i \in V}$;
and let us form the generating polynomial
\begin{equation}
   P_G({\bf t} , \boldsymbol{\lambda})
   \;=\;
   \sum_{A \subseteq E}  \left( \prod_{e \in A}  \lambda_e \right)
                         \left( \prod_{i \in V}  w_i(A) \right)
\end{equation}
where
\begin{alignat}{2}
   & \text{Lee--Yang:}
   &
   & w_i(A) \;=\;   {\begin{cases}
                        1    & \text{if $\deg_A(i)$ is even} \\
                        t_i  & \text{if $\deg_A(i)$ is odd}
                     \end{cases}
                    }
   \\[4mm]
   & \text{Heilmann--Lieb:}
   & \qquad
   & w_i(A) \;=\;   {\begin{cases}
                        1    & \text{if $\deg_A(i) = 0$} \\
                        t_i  & \text{if $\deg_A(i) = 1$} \\
                        0    & \text{if $\deg_A(i) > 1$}
                     \end{cases}
                    }
\end{alignat}
Then the Lee--Yang and Heilmann--Lieb theorems assert that
$P_G({\bf t} , \boldsymbol{\lambda}) \neq 0$
whenever $\lambda_e \ge 0$ for all $e$ and
$\real t_i > 0$ for all $i$ (or $\real t_i < 0$ for all $i$).
Indeed, the Lee--Yang and Heilmann--Lieb theorems can be given a unified
proof.\footnote{
   See \cite[Sections~4.7 and 4.8, and Remark 4 in Section~10.1]{Choe_HPP}.
   See also \cite{Ruelle_99a,Ruelle_99b} for related results,
   proven by a different method.
} 
The half-plane theorem for the spanning-tree polynomial
will be discussed in Section~\ref{sec.spanning} below.

For a prototypical example of a ``hard'' theorem,
consider the multivariate generating polynomial
for independent sets of vertices in a graph $G=(V,E)$:\footnote{
   We recall that a set $A \subseteq V$ is called
   \defword{independent}\/ (or \defword{stable}\/)
   if it does not contain any pair of adjacent vertices.
   This is totally unrelated to the concept of
   independent sets of elements in a matroid;
   it is unfortunate that the same word is used.
}
\begin{equation}
      Z_G(\bw)   \;=\;  \!\!\!\sum_{\begin{scarray}
                                      A \subseteq V \\
                                      A \, {\rm independent}
                                   \end{scarray}}
                       \prod\limits_{i \in A} w_i
  \label{eq.Zhardcore}
      \;,
\end{equation}
where $\bw = \{w_i\}_{i \in V}$ are (complex) vertex weights.
[In statistical mechanics, \reff{eq.Zhardcore} is called the
 grand partition function of the \defword{hard-core lattice gas}\/ on $G$;
 the weights $w_i$ are called \defword{activities}\/
 or \defword{fugacities}\/.]
It can then be shown \cite{Scott-Sokal} that if
$G$ is a graph of maximum degree $\Delta$,
and $|w_i| \le (\Delta-1)^{\Delta-1}/\Delta^\Delta$ for all $i \in V$,
then $Z_G(\bw) \neq 0$.
Indeed, this is a corollary of a more general result \cite{Scott-Sokal}
that provides a sufficient condition on a set of radii
${\bf R} = \{R_i\}_{i \in V}$
for $Z_G$ to be nonvanishing in the closed polydisc $|w_i| \le R_i$.\footnote{
   These results are, in turn, closely related to the
   Lov\'asz local lemma \cite{EL,ES,Spencer_75,Spencer_77}
   in probabilistic combinatorics;
   see \cite{Scott-Sokal} for a detailed development of this connection.
}

These theorems belong, in fact, to a long line of results
in the mathematical-physics literature asserting that
the grand partition function of a repulsive lattice gas
is nonvanishing in some polydisc at small fugacity.
(Physically this corresponds to the absence of phase transition
 for a gas at low density.)
Two main methods of proof are used:
\begin{itemize}
   \item[1)]  The coefficients in the Taylor expansion
      of $\log Z_G(\bw)$ are interpreted combinatorially:
      this is the so-called \defword{Mayer expansion}\/ \cite{Uhlenbeck_62}.
      Then, by direct estimation involving some nontrivial combinatorics,
      the Mayer expansion is proven to be convergent
      in a suitable polydisc $|w_i| < R_i$
      \cite{Penrose_67,Cammarota_82,Seiler_82,Brydges_86,Kotecky_86,%
Brydges_87,Brydges_88,Simon_93,%
Procacci_98,Procacci_99,Brydges_99,Miracle_00,Bovier_00}.
      It immediately follows that $Z_G$ is nonvanishing in this same polydisc.
   \item[2)]  The nonvanishing of $Z_G$ in the closed polydisc $|w_i| \le R_i$
      is proven directly, by induction on the number of vertices in $G$
      \cite{Dobrushin_96a,Dobrushin_96b,Sokal_chromatic_bounds,Scott-Sokal}.
      It then immediately follows that the Taylor series for $\log Z_G$
      is convergent in the interior of this polydisc.
\end{itemize}

The foregoing results have a wider application than one might
at first think, because the hard-core lattice gas
is not merely {\em one}\/ interesting statistical-mechanical model;
it is, in fact, the {\em universal}\/ statistical-mechanical model
in the sense that any statistical-mechanical model
living on a vertex set $V_0$ can be mapped onto a gas
of nonoverlapping ``polymers'' on $V_0$,
i.e.\ a hard-core lattice gas on the intersection graph of $V_0$
\cite[Section~5.7]{Simon_93}.\footnote{
   The {\em intersection graph}\/ of a finite set $S$
   is the graph whose vertices are the nonempty subsets of $S$,
   and whose edges are the pairs with nonempty intersection.
}
This construction, which is termed the
``polymer expansion'' or ``cluster expansion'',
is an important tool in mathematical statistical mechanics
\cite{Seiler_82,Brydges_86,Glimm_87,Brydges_99,Borgs_lectures};
it is widely employed to prove the absence of phase transition
at high temperature, low temperature, large magnetic field,
low density, or weak nonlinear coupling.
In the antiferromagnetic Potts model it can be used to prove
the nonvanishing of $Z_G(q,\bv)$ in a region $|q| > R$
(see Section~\ref{sec.largeq}).


\section{\mbox{\protect\boldmath $q=0$}:
             Spanning trees, electric circuits, and all that}
   \label{sec.spanning}

As explained in Section~\ref{subsec.q0},
the spanning-tree polynomial $T_G(\bv)$ can be obtained
from the multivariate Tutte polynomial $Z_G(q,\bv)$
via the double limit in which $q \to 0$ and $\bv$ is infinitesimal.
In this section we shall study some remarkable properties of
the spanning-tree polynomial and its generalization to matroids.

\subsection{The matrix-tree theorem}  \label{subsec.matrixtree}

Let $G=(V,E)$ be a connected graph,
and let $L_G(\bx)$ be the Laplacian matrix for $G$
with edge weights $\bx = \{x_e\}_{e \in E}$:
\begin{equation}
   L_G(\bx)
   \;=\;
   \begin{cases}
       - \sum_{e \sim ij} x_e                & \text{if $i \neq j$} \\
       \sum_{k \neq i} \sum_{e \sim ik} x_e  & \text{if $i=j$}
   \end{cases} 
 \label{def.Laplacian}
\end{equation}
where $e \sim ij$ denotes that $e$ connects $i$ to $j$.
[Equivalently, $L_G(\bx) = B X B^{\rm T}$
 where $B$ is the directed vertex-edge incidence matrix for
 any orientation of $G$, and $X = \diag(\bx)$.]
By construction, the row (and column) sums of $L_G(\bx)$ are all zero,
so its determinant is zero.
Now let $i$ be any vertex of $G$,
and let $L_G(\bx)_{\setminus i}$ be the matrix obtained from $L_G(\bx)$
by deleting the $i$th row and column.
Kirchhoff \cite{Kirchhoff_1847} proved in 1847
the following striking result
(see also \cite{Brooks_40,Nerode_61,Moon_70,Chen_71,Chaiken_78,Zeilberger_85}):

\begin{prop}[matrix-tree theorem]
  \label{prop.matrixtree}
%
$\; \det L_G(\bx)_{\setminus i}$ is independent of $i$ and\break
equals $T_G(\bx)$,
the generating polynomial of spanning trees in $G$.
\end{prop}

\noindent
Many different proofs of the matrix-tree theorem are now available;
one simple proof is based on the Cauchy--Binet theorem
in matrix theory (see e.g.\ \cite{Moon_70}).

More generally, it turns out that {\em each}\/ minor of $L_G(\bx)$
enumerates a suitable class of rooted spanning forests
\cite{Chaiken_82,Moon_94,Abdesselam_03}.
To formulate this result,
suppose that $I,J \subseteq V$,
and let us denote by $L_G(\bx)_{\setminus I,J}$
the matrix obtained from $L_G(\bx)$
by deleting the columns $I$ and the rows $J$;
when $I=J$, we write simply $L_G(\bx)_{\setminus I}$.
Then the ``principal-minors matrix-tree theorem'' states that
\begin{equation}
   \det L_G(\bx)_{\setminus \{i_1,\ldots,i_r\}}
   \;=\;
   \sum_{F \in {\cal F}(i_1,\ldots,i_r)} \, \prod_{e \in F} x_e   \;,
 \label{eq.principal}
\end{equation}
where the sum runs over all spanning forests $F$ in $G$
composed of $r$ disjoint trees, each of which contains exactly one
of the ``root'' vertices $i_1,\ldots,i_r$.
This theorem can easily be derived by applying
Proposition~\ref{prop.matrixtree} to the
graph in which the vertices $i_1,\ldots,i_r$
are contracted to a single vertex.
Furthermore, the ``all-minors matrix-tree theorem''
(whose proof is a bit more intricate,
 see \cite{Chaiken_82,Moon_94,Abdesselam_03})
states that for any subsets $I,J$ of the same cardinality $r$, we have
\begin{equation}
   \det L_G(\bx)_{\setminus I,J}  \;=\;
   \sum_{F \in {\cal F}(I|J)} \epsilon(F,I,J) \, \prod_{e \in F} x_e   \;,
 \label{eq.all-minors}
\end{equation}
where the sum runs over all spanning forests $F$ in $G$
composed of $r$ disjoint trees, each of which contains exactly one
vertex from $I$ and exactly one vertex (possibly the same one) from $J$;
here $\epsilon(F,I,J) = \pm 1$ are signs
whose precise definition is not needed here.

The virtue of the matrix-tree theorem is that
enumerative questions about spanning trees
(and {\em rooted}\/ spanning forests)
can be reduced to linear algebra.

\subsection{Electric circuits and the half-plane property for graphs}

Now let us consider the graph $G$ as an electrical network:
to each edge $e$ we associate a complex number $x_e$,
called its {\em conductance}\/ (or {\em admittance}\/).\footnote{
   Complex admittances are relevant, for instance,
   in alternating-current circuits.
   Indeed, the reader is invited to imagine that we are considering
   an alternating-current circuit at some fixed frequency $\omega$.
}
[The \defword{resistance}\/ (or \defword{impedance}\/) is $1/x_e$.]
Suppose that we inject currents ${\bf J} = \{ J_i \} _{i \in V}$
into the vertices.
What node voltages $\bvarphi = \{ \varphi_i \} _{i \in V}$
will be produced?
Applying Kirchhoff's law of current conservation at each vertex
and Ohm's law on each edge, it is not hard to see that
the node voltages and current inflows satisfy the linear system
\begin{equation}
   L_G(\bx) \, \bvarphi  \;=\;  {\bf J}
   \;.
 \label{eq.electric}
\end{equation}
It is then natural to ask:
Under what conditions does this system have a (unique) solution?
Two obvious constraints arise from the fact that
the row and column sums of $L_G(\bx)$ are zero:
firstly, the current vector must satisfy $\sum_{i \in V} J_i = 0$
(``conservation of total current''), or else no solution will exist;
and secondly, if $\bvarphi$ is any solution,
then so is $\bvarphi + c{\bf 1}$ for any $c$
(``only voltage {\em differences}\/ are physically observable'').
So let us assume that $\sum_{i \in V} J_i = 0$;
and let us break the redundancy in the solution
by fixing the voltage to be zero
at some chosen reference node $i_0 \in V$ (``ground'').
Does the modified system
\begin{equation}
   L_G(\bx)_{\setminus i_0} \, \bvarphi_{\setminus i_0}
  \;=\;  {\bf J}_{\setminus i_0}
 \label{eq.electric2}
\end{equation}
then have a unique solution?
This will be so if and only if $\det L_G(\bx)_{\setminus i_0}$
is nonzero --- which, by virtue of the matrix-tree theorem,
is equivalent to $T_G(\bx)$ being nonzero.

Simple counterexamples show that this is not always the case.
Suppose, for instance, that $G$ consists of a pair of vertices
connected by two edges $e,f$ in parallel.
With $x_e = 1$ and $x_f = -1$, 
it is easy to see that no solution exists (except when ${\bf J} = 0$).
Of course, the reader will object that negative resistances
are physically unrealizable!
Fine:  consider instead $x_e = i$ and $x_f = -i$;
once again there is no solution (unless ${\bf J} = 0$).
This example corresponds physically to a perfectly lossless capacitor
together with a perfectly lossless inductor,
exactly at their resonant frequency.
But once again, the reader will rightly object that {\em perfectly}\/
lossless components are physically unrealizable;
every component in the real world exhibits {\em some}\/ dissipation.
This reasoning leads us to conjecture, on physical grounds,
that if $\real x_e > 0$ for all $e$
(i.e.\ every branch is strictly dissipative),
then the network is uniquely solvable once we fix the voltage
at a single reference node $i_0 \in V$.
This conjecture turns out to be true, and we have:

\begin{theorem}
 \label{thm1.1}
Let $G$ be a connected graph.
Then the spanning-tree polynomial $T_G$ has the ``half-plane property'' (HPP),
i.e.\ $\real x_e > 0$ for all $e$ implies $T_G(\bx) \neq 0$.
\end{theorem}

The proof of Theorem~\ref{thm1.1} is not difficult:
Consider any nonzero complex vector
$\bvarphi = \{ \varphi_i \} _{i \in V}$ satisfying $\varphi_{i_0} = 0$.
Because $G$ is connected, we have $B^{\rm T} \bvarphi \neq 0$.
Therefore, the quantity
\begin{equation}
   \bvarphi^* L_G(\bx) \bvarphi  \;=\;
   \bvarphi^* B X B^{\rm T} \bvarphi  \;=\;
   \sum_{e \in E} |(B^{\rm T} \bvarphi)_e|^2 \, x_e
 \label{energy_form}
\end{equation}
has strictly positive real part whenever $\real x_e > 0$ for all $e$;
so in particular $(B X B^{\rm T} \bvarphi)_i \neq 0$ for some $i \neq i_0$.
It follows that the submatrix of $L_G(\bx)$ obtained by suppressing
the $i_0$th row and column is nonsingular,
and so has a nonzero determinant.
Theorem~\ref{thm1.1} then follows from the matrix-tree theorem.\footnote{
   This proof is well known in the circuit-theory literature:
   see e.g.\ \cite[Section 2.7]{Chen_71}
   as well as the related results in
   \cite[pp.~398--401, 430--431 and 850--851]{Desoer_69}
   \cite[pp.~52--53 and 67--69]{Penfield_70}.
   It has, moreover, a natural physical interpretation:
   if $\bvarphi = \{ \varphi_i \} _{i \in V}$ are the node voltages,
   then the real part of the quadratic form \reff{energy_form}
   is the total power dissipated in the circuit.
}

It cannot be overemphasized how remarkable this theorem is.
Suppose, for instance, that $G$ is the complete graph on $n$ vertices;
then $T_G(\bx)$ is a homogeneous polynomial of degree $n-1$
in $n(n-1)/2$ variables, containing $n^{n-2}$ monomials.
It seems utterly miraculous, at first sight,
that a polynomial of such complexity should be nonvanishing
whenever $\real x_e > 0$ for all $e$.
The fact that it is so clearly expresses a deep property
arising from the underlying combinatorial structure.

An immediate corollary of Theorem~\ref{thm1.1}
is that the complementary spanning-tree polynomial
\begin{equation}
   \widetilde{T}_G(\bx)  \;=\;
   \left( \prod\limits_{e \in E} x_e \right) \, T_G(1/\bx)
\end{equation}
also has the half-plane property,
since the map $x_e \mapsto 1/x_e$ takes the right half-plane onto itself.

\subsection{The half-plane property for matroids}

{}From a combinatorial point of view,
the noteworthy fact is that the spanning trees of $G$
constitute the bases of the graphic matroid $M(G)$,
and their complements constitute the bases of the cographic matroid $M^*(G)$.
So $T_G$ and $\widetilde{T}_G$ are the
multivariate basis generating polynomials for $M(G)$ and $M^*(G)$,
respectively.
This naturally suggests generalizing Theorem~\ref{thm1.1}
to more general matroids and, perhaps, to more general set systems.
Before posing these questions precisely,
we need to fix some notation and terminology.

A \defword{set system}\/ (or \defword{hypergraph}\/) $\scrs$
on the (finite) ground set $E$
is simply a collection $\scrs$ of subsets of $E$.
Given any set system $\scrs$ on $E$,
we define its {\em (multivariate) generating polynomial}\/ to be
\begin{equation}
   P_\scrs(x)   \;=\;  \sum_{S \in \scrs} \, \prod_{e \in S} x_e  \;.
\end{equation}
The {\em rank}\/ of a set system is the maximum cardinality
of its members (by convention we set rank $= -\infty$ if $\scrs = \emptyset$);
equivalently, it is the degree of the generating polynomial $P_\scrs$.
A set system $\scrs$ is {\em $r$-uniform}\/ if $|S|=r$ for all $S \in \scrs$,
or equivalently if its generating polynomial $P_\scrs$ is
homogeneous of degree $r$.
We shall be particularly interested in the
\defword{basis generating polynomial}\/ of a matroid $M$,
\begin{equation}
   P_{\scrb(M)}(x)  \;=\; B_M(x)
       \;=\;  \sum_{B \in \scrb(M)}  \, \prod_{e \in B}  x_e
   \;.
\end{equation}

We can now pose the following questions concerning
possible extensions of Theorem~\ref{thm1.1}:

\begin{qn}
  \label{question1.2}
For which matroids $M$ does the basis generating polynomial $P_{\scrb(M)}$
have the half-plane property?
\end{qn}

\noindent
More generally:

\begin{qn}
  \label{question1.3}
For which $r$-uniform set systems $\scrs$ does the generating polynomial
$P_{\scrs}$ have the half-plane property?
\end{qn}

These questions were recently studied in a long paper by
Choe, Oxley, Wagner and myself \cite{Choe_HPP}.
Our original conjecture was that all matroids
(and no non-matroidal set systems) have the half-plane property.
That would be nice and neat, but it turns out to be false;
and the truth is considerably more interesting and subtle.
Our conjecture is half right:
an $r$-uniform set system with the half-plane property
is necessarily the set of bases of a matroid \cite[Theorem~7.1]{Choe_HPP}.
But not every matroid has the half-plane property,
and we do not yet have a complete characterization of those that do.
Nevertheless, we can find large classes of matroids
with the half-plane property:
\begin{itemize}
   \item[(a)]  Every complex unimodular matroid
--- or what is equivalent, sixth-root-of-unity matroid \cite{Whittle_97} ---
has the half-plane property
\cite[Theorem~8.1 and Corollary~8.2]{Choe_HPP}.\footnote{
    It is proven in \cite[Theorem~8.9]{Choe_HPP}
    that a matroid is complex unimodular if and only if
    it is sixth-root-of-unity.
}
This class properly includes the regular matroids,
which in turn properly includes the graphic and cographic matroids.
The proof of the half-plane property for complex unimodular matroids
is, in fact, a direct generalization of the proof just given
for Theorem~\ref{thm1.1}.
   \item[(b)]  Every uniform matroid has the half-plane property.
(Indeed, every loopless uniform matroid has the Brown--Colbourn property,
 which is stronger than the half-plane property:
 see Section~\ref{sec.reliability} below.)
   \item[(c)]  A significant subclass of transversal matroids
--- those we call ``nice'' ---
have the half-plane property \cite[Section~10]{Choe_HPP}.\footnote{
   In \cite[Conjecture 13.16 and Question 13.17]{Choe_HPP}
   it was conjectured that all rank-3 transversal matroids have
   the half-plane property, and the question was raised whether
   it might even be true that {\em all}\/ transversal matroids
   have the half-plane property.
   Both conjectures are in fact false.
   Choe and Wagner \cite{Choe_Rayleigh} have exhibited a
   rank-4 transversal matroid that lacks the half-plane property
   (and indeed lacks the weaker Rayleigh property,
    see Section~\ref{subsec.Rayleigh} below).
   I~have found a class of rank-3 transversal matroids that
   lack the half-plane property,
   and plan to publish them elsewhere \cite{Sokal_hpp2}.
}
   \item[(d)]  All matroids of rank or corank at most 2 have the
half-plane property \cite[Corollary~5.5]{Choe_HPP},
as do all matroids on a ground set of at most 6 elements
\cite[Proposition~10.4]{Choe_HPP}.
   \item[(e)]  The class of matroids with the half-plane property
is closed under minors, duality, direct sums, 2-sums,
series and parallel connection, full-rank matroid union,
and some special cases of principal truncation, principal extension,
principal cotruncation and principal coextension \cite[Section~4]{Choe_HPP}.
\end{itemize}
Furthermore, we show that:
\begin{itemize}
   \item[(f)]  A binary matroid has the half-plane property
       if and only if it is regular
       \cite[Corollary~8.16]{Choe_HPP}.
   \item[(g)]  A ternary matroid has the half-plane property
       if and only if it is a sixth-root-of-unity matroid
       \cite[Corollary~8.17]{Choe_HPP}.
\end{itemize}
Finally, we can show that certain matroids do {\em not}\/ have the
half-plane property:  among these are the Fano matroid $F_7$,
the non-Fano matroid $F_7^-$,
their relaxations $F_7^{--}$, $F_7^{-3}$ and $M(K_4)+e$,
the matroids $P_8$, $P'_8$ and $P''_8$,
the Pappus and non-Pappus matroids,
the free extension (non-Pappus $\drop\,9)+e$,
and all their duals
\cite[Section~11]{Choe_HPP}.
The first six of these examples are minor-minimal,
and we conjecture that the others are as well;
but we strongly suspect that this list is incomplete,
and indeed we consider it likely that the set of minor-minimal
non-half-plane-property matroids is infinite.

More generally, we consider homogeneous multiaffine polynomials
\begin{equation}
   P(\bx)  \;=\; \sum_{S \subseteq E, |S| = r} \! a_S \, \prod_{e \in S} x_e
\end{equation}
with arbitrary complex coefficients $a_S$ (not necessarily 0 or 1).
We prove two {\em necessary}\/ conditions for $P \not\equiv 0$
to have the half-plane property:
\begin{itemize}
   \item[(a)] $P$ must have the ``same-phase property'',
      i.e.\ all the nonzero coefficients $a_S$ must have the same phase
      \cite[Theorem~6.1]{Choe_HPP}.
      So without loss of generality we can assume that all the $a_S$
      are real and nonnegative.
   \item[(b)] The {\em support}\/
      $\supp(P) = \{S \subseteq E \colon\;  a_S \neq 0 \}$
      must be the collection of bases of a matroid
      \cite[Theorem~7.1]{Choe_HPP}.
\end{itemize}
This latter fact is particularly striking:  it shows that matroids
arise {\em naturally}\/ from a consideration of homogeneous multiaffine
polynomials with the half-plane property.
We do not know whether the converse of this result
is true, i.e.\ whether for every matroid $M$ there exists a
homogeneous multiaffine polynomial $P$ with the half-plane property
such that $\supp(P) = \scrb(M)$.
But it is true, at least, for all matroids representable over $\C$
\cite[Corollary~8.2]{Choe_HPP}.

We also also give two {\em sufficient}\/ conditions
for a homogeneous multiaffine polynomial $P$
to have the half-plane property (or be identically zero):
\begin{itemize}
   \item[(a)] {\em Determinant condition}\/ \cite[Theorem~8.1]{Choe_HPP}:
      $a_S = |\det(A \restrict S)|^2$ for some $r \times n$ complex matrix $A$
      [here $n=|E|$, and $A \restrict S$ denotes the square submatrix of $A$
       using the columns indexed by the set $S$].
      This corresponds to $P(\bx) = \det(A X A^*)$
      where $X = \diag(\bx)$ and ${}^*$ denotes Hermitian conjugate.
   \item[(b)] {\em Permanent condition}\/ \cite[Theorem~10.2]{Choe_HPP}:
      $a_S = \per(\Lambda \restrict S)$ for some $r \times n$ nonnegative
      matrix $\Lambda$.
      This corresponds to $P(\bx) = \per(\Lambda X)$.
\end{itemize}
Unfortunately, the relationship between these sufficient conditions
and the half-plane property looks complicated.
Neither family of polynomials contains the other;
their intersection is nonempty;
and their union is a proper subset of the set of
all homogeneous multiaffine polynomials with the half-plane property.

In any case, the matroids with the half-plane property
form a very natural class, which deserves further study.
A long list of open questions can be found in \cite[Section~13]{Choe_HPP}.

\subsection{The Rayleigh property for graphs and matroids}
   \label{subsec.Rayleigh}

Consider once again the graph $G$ as an electrical network,
with conductances $\{x_e\}_{e \in E}$ on the edges,
and select a pair of distinct vertices $i,j \in V$.
Then Kirchhoff \cite{Kirchhoff_1847} showed in 1847 that
the effective conductance between nodes $i$ and $j$ is
\begin{equation}
   \scry_{ij}({\bf x})  \;=\;   {T_G(\bx) \over T_{G/ij}(\bx)}
   \;,
  \label{eq.effcond}
\end{equation}
where $G/ij$ is the graph obtained from $G$ by contracting
$i$ and $j$ to a single vertex.\footnote{
   {\sc Proof:}  By definition, $1/\scry_{ij}({\bf x})$ is the
   voltage difference induced between nodes $i$ and $j$
   if we inject 1 ampere of current into $i$
   and extract 1 ampere from $j$.
   So let node $j$ be ``ground'', and apply \reff{eq.electric2}
   with $i_0 = j$ and ${\bf J} = \delta_i - \delta_j$.  We then have
   $$
      1/\scry_{ij}({\bf x})  \;=\;
         [L_G(\bx)_{\setminus j}]^{-1}_{ii} \\[2mm]
      \;=\;
         {\det L_G(\bx)_{\setminus \{i,j\}}
          \over
          \det L_G(\bx)_{\setminus j}
         }
   $$
   by Cramer's rule.
   The principal-minors matrix-tree theorem \reff{eq.principal}
   then tells us that
   $\det L_G(\bx)_{\setminus j} = T_G(\bx)$ and
   $\det L_G(\bx)_{\setminus \{i,j\}} = T_{G/ij}(\bx)$.
}

Suppose now (in contrast to the preceding subsections)
that all the conductances $x_e$ are real and positive.
Physical intuition suggests that increasing one or more
of the branch conductances $x_e$ cannot cause the effective conductance
between any pair of nodes to {\em decrease}\/;
in other words, we conjecture that $\scry_{ij}(\bx)$ is a nondecreasing
function of all the $x_e$ on $[0,\infty)^E$.
Let us rephrase this statement by considering the graph $H$
obtained from $G$ by adjoining a new edge $f$ connecting $i$ to $j$;
then $G = H \setminus f$ and $G/ij = H/f$.  Using
\begin{subeqnarray}
   T_{H \setminus f}(\bx)  & = &
     T_{H \setminus \{e,f\}}(\bx) \,+\, x_e T_{(H \setminus f)/e}(\bx)  \\
   T_{H/f}(\bx)  & = &
     T_{(H/f) \setminus e}(\bx) \,+\, x_e T_{H / \{e,f\}}(\bx)
\end{subeqnarray}
we find
\begin{subeqnarray}
   {\partial \over \partial x_e} \, \scry_{ij}({\bf x})
   & = &
   { T_{(H/e) \setminus f}(\bx) \, T_{(H/f) \setminus e}(\bx) \;-\;
     T_{H / \{e,f\}}(\bx) \, T_{H \setminus \{e,f\}}(\bx)
        \over
     T_{H/f}(\bx)^2
   }
    \qquad \\[2mm]
   & = &
   { T_{H/e}(\bx) \, T_{H/f}(\bx) \;-\;
     T_{H / \{e,f\}}(\bx) \, T_H(\bx)
        \over
     T_{H/f}(\bx)^2
   }
   \;\,.
\end{subeqnarray}
The conjecture then becomes:

\begin{theorem}[Rayleigh property for graphs]
 \label{thm.Rayleigh_graphs}
Let $H=(V,E)$ be a finite undirected graph,
and let $e,f \in E$ with $e \neq f$.
Then
\begin{equation}
   T_{H/e}(\bx) \, T_{H/f}(\bx) \;-\; T_{H / \{e,f\}}(\bx) \, T_H(\bx)
   \;\ge\; 0
 \label{eq.rayleigh.graph}
\end{equation}
for all $\bx \ge 0$.
\end{theorem}

\noindent
Please note that \reff{eq.rayleigh.graph} has a beautiful
probabilistic interpretation:
Consider the probability measure on spanning trees $T \subseteq H$
giving weight $\big( \prod_{e \in T} x_e \big) \big/ T_H(\bx)$
to the tree $T$.
Then \reff{eq.rayleigh.graph} asserts that the events
$e \in T$ and $f \in T$ are {\em negatively correlated}\/.

Theorem~\ref{thm.Rayleigh_graphs} has both algebraic and
variational/probabilistic proofs:
see e.g.\ \cite[Section~3.8]{Balabanian_69}
\cite[Theorem~2.1]{Feder-Mihail}
\cite[Theorems~5.2 and 5.6]{Choe_Rayleigh}
for the former
and \cite{Doyle-Snell,Pemantle_95,Lyons_98,Lyons-Peres} for the latter.

In an important recent paper, Choe and Wagner \cite{Choe_Rayleigh}
have investigated the generalization of this problem to matroids.
Let $M$ be a matroid with ground set $E$,
and let $B_M(\bx) = P_{\scrb(M)}(\bx)$
be its (multivariate) basis generating polynomial.
We say that $M$ is a \defword{Rayleigh matroid}\/
(or has the \defword{Rayleigh property}\/)
in case
\begin{equation}
   B_{M/e}(\bx) \, B_{M/f}(\bx) \;-\; B_{M / \{e,f\}}(\bx) \, B_M(\bx)
   \;\ge\; 0
 \label{eq.rayleigh.matroid}
\end{equation}
for all $e,f \in E$ with $e \neq f$
and all $\bx \ge 0$.
Choe and Wagner \cite{Choe_Rayleigh} have shown that:
\begin{itemize}
   \item[(a)]  The class of Rayleigh matroids
      is closed under minors, duality, direct sums, and 2-sums.
   \item[(b)]  A binary matroid is Rayleigh if and only if
      it does not contain $S_8$ as a minor.
      Equivalently, a binary matroid is Rayleigh if and only if
      it can be constructed from regular matroids, $F_7$, $F_7^*$
      and $AG(3,2)$ by taking direct sums and 2-sums.
   \item[(c)]  Every matroid with the half-plane property
       is Rayleigh (but the converse is false).
   \item[(d)]  Every matroid on a ground set of at most 7 elements
       is Rayleigh.
\end{itemize}
More recently, Wagner \cite{Wagner_rank3} has also shown that
\begin{itemize}
   \item[(e)]  Every matroid of rank or corank at most 3
       is Rayleigh.
\end{itemize}
On the other hand, Choe and Wagner \cite{Choe_Rayleigh}
have exhibited a rank-4 transversal matroid that is not Rayleigh.

The most striking of these results, in my opinion, are (b) and (c).
For binary matroids, (b) provides a pair of beautiful
characterizations of the Rayleigh property:
one by excluded minors and the other by internal structure.
It is at least conceivable that analogous results
could be obtained for ternary matroids,
$GF(4)$-representable matroids,
or even arbitrary matroids.
As for (c), one might have thought, {\em a priori}\/,
that the half-plane property and the Rayleigh property
are distinct but unrelated properties of electric circuits
(one for complex edge conductances, the other for positive real
edge conductances), which thus hold for all graphic matroids
but extend to distinct and unrelated families of non-graphic matroids.
Item (c) shows that the truth is quite different:
the half-plane property is strictly stronger than the Rayleigh property.
Indeed, in some vague sense it seems that the half-plane property
is {\em quite a bit}\/ stronger than the Rayleigh property:
non-HPP matroids seem to be plentiful,
while non-Rayleigh matroids seem to be fairly rare.
In any case, the Rayleigh matroids form a very natural class
and deserve further study;
see \cite{Choe_Rayleigh} for a list of open problems.

\section{\mbox{\protect\boldmath $q=0$} again:  The reliability polynomial}
  \label{sec.reliability}

\subsection{The Brown--Colbourn property for graphs}

Once again let $G=(V,E)$ be a connected graph,
considered now as a communications network
with unreliable communication channels:
the edge $e$ is assumed to be operational with probability $p_e$
and failed with probability $1-p_e$,
independently for each edge.
Let $R_G(\bp)$ be the probability that
every node is capable of communicating with every other node
(this is the so-called \defword{all-terminal reliability}\/).
Clearly we have
\begin{equation}
   R_G(\bp)  \;=\;
   \sum_{\begin{scarray}
           A \subseteq E \\
           (V,A) \, {\rm connected}
         \end{scarray}}
   \prod_{e \in A} p_e  \prod_{e \in E \setminus A} (1-p_e)
   \;,
\end{equation}
where the sum runs over all connected spanning subgraphs of $G$,
and we have written $\bp = \{p_e\}_{e \in E}$.
We call $R_G(\bp)$ the (multivariate) {\em reliability polynomial}\/
\cite{Colbourn_87} for the graph $G$;
it is a multiaffine polynomial,
i.e.\ of degree at most 1 in each variable separately.
If the edge probabilities $p_e$ are all set to the same value $p$,
we write the corresponding univariate polynomial as $R_G(p)$,
and call it the univariate reliability polynomial.
We are interested in studying the zeros of these polynomials
when the variables $p_e$ (or $p$) are taken to be {\em complex}\/ numbers.

The reliability polynomial $R_G(\bp)$ is, of course, just the
connected-spanning-subgraph polynomial $C_G(\bv)$ in disguise:
\begin{eqnarray}
   R_G(\bp)  & = &   \left[ \prod_{e \in E} (1-p_e) \right]
                      C_G \!\left( { \bp \over {\bf 1} - \bp} \right)
        \\[4mm]
   C_G(\bv)  & = &   \left[ \prod_{e \in E} (1+v_e) \right]
                      R_G \!\left( { \bv \over {\bf 1} + \bv} \right)
\end{eqnarray}
So we can equally well work with $C_G(\bv)$,
making the change of variables $v_e = p_e/(1-p_e)$.

Brown and Colbourn \cite{Brown_92}
studied a number of examples and made the following conjecture:
\begin{quote}
   {\bf Univariate Brown--Colbourn conjecture.}
   For any connected graph $G$,
   the zeros of the univariate reliability polynomial $R_G(p)$
   all lie in the closed disc $|p-1| \le 1$.
   In other words, if $|p-1| > 1$, then $R_G(p) \neq 0$.
\end{quote}
Subsequently, I~proposed \cite{Sokal_chromatic_bounds}
a multivariate extension of the Brown--Colbourn conjecture:
\begin{quote}
   {\bf Multivariate Brown--Colbourn conjecture.}
   For any connected graph $G$,
   if $\mbox{$|p_e-1| > 1$}$ for all edges $e$, then $R_G(\bp) \neq 0$.
\end{quote}
In terms of $C_G$, the multivariate Brown--Colbourn conjecture states that
if $G$ is a {\em loopless}\/ connected graph
and $|1+v_e| < 1$ for all edges $e$, then $C_G(\bv) \neq 0$.
Loops must be excluded because a loop $e$ multiplies $C_G$
by a factor $1+v_e$ but leaves $R_G$ unaffected.

A few years ago, Wagner \cite{Wagner_00} proved,
using an ingenious and complicated construction,
that the univariate Brown--Colbourn conjecture
holds for all series-parallel graphs.\footnote{
   Unfortunately, there seems to be no completely standard
   definition of ``series-parallel graph'';
   a plethora of slightly different definitions can be found in the literature
   \cite{Duffin_65,Colbourn_87,Oxley_86,Oxley_92,Brandstadt_99}.
   So let me be completely precise about my own usage:
   I shall call a loopless graph {\em series-parallel}\/
   if it can be obtained from a forest by a finite sequence of
   series and parallel extensions of edges
   (i.e.\ replacing an edge by two edges in series or two edges in parallel).
   I shall call a general graph (allowing loops) series-parallel
   if its underlying loopless graph is series-parallel.
   Some authors write ``obtained from a tree'', ``obtained from $K_2$''
   or ``obtained from $C_2$'' in place of ``obtained from a forest'';
   in my terminology these definitions yield, respectively,
   all {\em connected}\/ series-parallel graphs,
   all connected series-parallel graphs whose blocks form a path,
   or all {\em 2-connected}\/ series-parallel graphs.
   See \cite[Section 11.2]{Brandstadt_99} for a more extensive bibliography.
}
As an illustration of the power of the ``multivariate ideology'',
let us now prove the stronger result \cite{Sokal_chromatic_bounds}
that the multivariate Brown--Colbourn conjecture
holds for all series-parallel graphs ---
and let us moreover do it in two lines:\footnote{
   This proof can be found in
   \cite[Remark 3 in Section 4.1]{Sokal_chromatic_bounds}
   or in \cite[Theorem 5.6(c) $\Longrightarrow$ (a)]{Royle-Sokal}.
}

\begin{theorem}
Let $G$ be a loopless connected series-parallel graph.
If $|1+v_e| < 1$ for all edges $e$, then $C_G(\bv) \neq 0$.
\end{theorem}

\begin{proof}
How should we prove a theorem for series-parallel graphs?
The answer is obvious: just show that series and parallel reduction
preserve the hypothesis $|1+v_e| < 1$.
If edges $e_1,e_2$ (with weights $v_1,v_2$) are in parallel,
they can be replaced by a single edge with weight $v_*$
satisfying $1+v_* = (1+v_1)(1+v_2)$;
obviously $|1 + v_i| < 1$ for $i=1,2$ implies that $|1 + v_*| < 1$.
Likewise, if edges $e_1,e_2$ (with weights $v_1,v_2$) are in series,
they can be replaced by a single edge with weight $v_*$
satisfying $1/v_* = 1/v_1 + 1/v_2$,
provided that we multiply $C_G$ by the prefactor $v_1 + v_2$
[cf.\ \reff{eq.series1}/\reff{eq.series2} specialized to $q=0$].
Now $|1+v| < 1$ is equivalent to $\real(1/v) < -1/2$,
so $\real(1/v_i) < -1/2$ for $i=1,2$ implies that $\real(1/v_*) < -1 < -1/2$;
furthermore, the prefactor $v_1 + v_2$ in \reff{eq.series2} is nonzero.
Since the multivariate Brown--Colbourn conjecture manifestly
holds for the base case of trees,
it necessarily holds for all loopless connected series-parallel graphs.
\end{proof}

Does the multivariate Brown--Colbourn conjecture hold for all graphs?
For several years I~would have bet that it does
(though I~was unable to find a proof);
but in 2002, Gordon Royle sent me the surprising news that
the multivariate Brown--Colbourn conjecture fails already
for the simplest non-series-parallel graph,
namely the complete graph $K_4$.
The construction turns out to be very simple \cite{Royle-Sokal}.
Since the {\em uni}\/variate Brown--Colbourn conjecture holds for $K_4$,
let us try the next simplest situation, namely the {\em bi}\/variate one
in which the six edges receive two different weights $a$ and $b$.
There are five cases:
\begin{itemize}
   \item[(a)]  One edge receives weight $a$ and the other five
      receive weight $b$:
\begin{equation}
   C_{K_4}(a,b)  \;=\;  (8b^3 + 5b^4 + b^5) + (8b^2 + 10b^3 + 5b^4 + b^5) a
\end{equation}
   \item[(b)]  A pair of nonintersecting edges receive weight $a$
      and the other four edges receive weight $b$:
\begin{equation}
   C_{K_4}(a,b)  \;=\;  (4b^3 + b^4) + (8b^2 + 8 b^3 + 2b^4) a
                           + (4b + 6b^2 + 4 b^3 + b^4) a^2
\end{equation}
   \item[(c)]  A pair of intersecting edges receive weight $a$
      and the other four edges receive weight $b$:
\begin{equation}
   C_{K_4}(a,b)  \;=\;  (3b^3 + b^4) + (10b^2 + 8 b^3 + 2b^4) a
                           + (3b + 6b^2 + 4 b^3 + b^4) a^2
\end{equation}
   \item[(d)]  A 3-star receives weight $a$ and the complementary
      triangle receives weight $b$:
\begin{equation}
   C_{K_4}(a,b)  \;=\;  (9b^2 + 3b^3) a + (6b + 9b^2 + 3b^3) a^2
                           + (1 + 3b + 3b^2 + b^3) a^3
\end{equation}
   \item[(e)]  A three-edge path receives weight $a$ and the
      complementary three-edge path receives weight $b$:
\begin{equation}
   C_{K_4}(a,b)  \;=\;  b^3 + (7b^2+3b^3)a + (7b+9b^2+3b^3)a^2
                           + (1+3b+3b^2+b^3)a^3
\end{equation}
\end{itemize}
Now plot the roots $a$ when $b$ traces out the circle $|1+b|=1$,
or vice versa.
In cases (b) and (d) it turns out that the roots
can enter the ``forbidden discs'' $|1+a| < 1$ and $|1+b| < 1$
(see \cite[Figures 1 and 2]{Royle-Sokal}).
In fact, this is quite easy to understand analytically:
if we solve the equation $C_{K_4}(a,b) = 0$ for $a$ in terms of $b$,
expanding in power series for $b$ near 0,
in cases (b) and (d) we obtain
\begin{equation}
   a   \;=\;  \delta_1 b + \delta_2 b^{3/2} + O(b^2)
\end{equation}
with $\delta_1 < 0$ and $\delta_2 \neq 0$.
If we now set $b = -1 + e^{i\theta}$,
then one of the roots will have $|1+a| < 1$ for small $\theta \neq 0$
(arising from the $\delta_2 b^{3/2}$ term).
A small perturbation (so that $|1+b| < 1$)
then yields a counterexample to the multivariate Brown--Colbourn conjecture.

In fact, with a little more work Royle and I are able to prove
the following \cite[Theorem 5.6]{Royle-Sokal}:

\begin{theorem}
A loopless connected graph $G$ has the multivariate Brown--Colbourn
property {\em if and only if} it is series-parallel.
\end{theorem}

Moreover, as a corollary of the failure of the multivariate Brown--Colbourn
conjecture for the complete graph $K_4$,
Royle and I are able to show \cite{Royle-Sokal}
that the univariate conjecture is false as well:
counterexamples include
a 4-vertex, 16-edge planar graph that can be obtained from $K_4$
by adding parallel edges, and a 1512-vertex, 3016-edge simple planar graph
that can be obtained from $K_4$ by adding parallel edges and then
subdividing edges.
These counterexamples to the univariate conjecture are fairly easy to find
once one has in hand the $K_4$ counterexample to the multivariate conjecture,
but would probably not otherwise be guessed.
This illustrates once again the advantages of considering
the multivariate problem, even if one is ultimately interested
in a particular univariate specialization.

\subsection{The Brown--Colbourn property for matroids}

The foregoing considerations can be extended to matroids.
Let $M$ be a matroid with ground set $E$,
and let $S_M(\bv)$ be its spanning-set polynomial.
We then say that $M$ has the
(multivariate) \defword{Brown--Colbourn property}\/
in case $S_M(\bv) \neq 0$ whenever $|1+v_e| < 1$ for all $e \in E$.\footnote{
   Unfortunately, in \cite{Choe_HPP} we used the dual definition
   (i.e., a matroid $M$ has the Brown--Colbourn property
    in the sense of \cite{Choe_HPP}
    if and only if $M^*$ has the Brown--Colbourn property as defined here).
   I~think the present definition is more natural
   in view of the relation to graphs
   (namely, $G$ has the Brown--Colbourn property as defined here
    if and only if its cycle matroid $M(G)$ does).
   I~therefore pledge to stick to this definition in the future.
}

Gordon Royle and I are currently investigating the Brown--Colbourn property
for matroids \cite{Royle-Sokal_matroids}.
Here are some of our results so far:
\begin{itemize}
   \item[(a)]  The class of matroids with the Brown--Colbourn property
      is closed under deletion,
      direct sums, 2-sums, series connection and parallel connection.\footnote{
          The class of matroids with the B--C property
          fails to be closed under contraction, at least for a trivial reason:
          namely, contraction of one of a set of parallel elements
          will produce one or more loops, which are incompatible
          with the B--C property [see item (c) below].
          But we do not know whether the contraction of a {\em simple}\/
          matroid with the B--C property can fail to have the B--C property.
}
   \item[(b)]  If $M$ has the Brown--Colbourn property,
      then it has the half-plane property as well \cite{Choe_HPP}.
      But the converse does not hold, as the example of $M(K_4)$ shows.
   \item[(c)]  If $M$ has the Brown--Colbourn property,
      it must be loopless.
   \item[(d)]  Every loopless uniform matroid
      (i.e., every uniform matroid $U_{r,n}$ with $1 \le r \le n$)
      has the Brown--Colbourn property \cite{Wagner_00,Choe_HPP}.
   \item[(e)]  A binary matroid has the Brown--Colbourn property
      if and only if it is the cycle matroid
      of a loopless series-parallel graph.
\end{itemize}
Most of these results are fairly easy to prove.

\section{\mbox{\protect\boldmath $q=0$} yet again:  Spanning forests}

As discussed in Section~\ref{subsec.matrixtree},
Kirchhoff's matrix-tree theorem
\cite{Kirchhoff_1847,Brooks_40,Nerode_61,Moon_70,Chaiken_78,Zeilberger_85}
and its generalizations to arbitrary minors
\cite{Chaiken_82,Moon_94,Abdesselam_03}
express the generating polynomials of spanning trees
and rooted spanning forests in a graph $G$
as determinants associated to the graph's Laplacian matrix.
Here I~would like to discuss a recently discovered extension of the
matrix-tree theorem that provides a compact representation
also for the generating polynomial of {\em unrooted}\/ spanning forests,
$F_G(\bw)$.

Recall first that one useful formalism for manipulating determinants
is Grassmann algebra;
indeed, any determinant can be written as a Gaussian ``integral''
over Grassmann variables.\footnote{
   For introductions to Grassmann algebra and Grassmann--Berezin integration,
   see e.g.\ \cite{Zinn-Justin,Abdesselam_03}.
}
In particular, it follows from the matrix-tree theorem that the
generating polynomials of spanning trees and rooted spanning forests
can be written as Gaussian Grassmann integrals.

Very recently, Caracciolo, Jacobsen, Saleur, Sportiello and I
\cite{Caracciolo_matrixtree}
have proven a generalization of the matrix-tree theorem
in which a large class of combinatorial objects are represented
by suitable {\em non-Gaussian}\/ Grassmann integrals.
As a special case, we show that the generating polynomial of
unrooted spanning forests,
which arises as a $q \to 0$ limit of the multivariate Tutte polynomial
[cf.\ \reff{scheme.qto0}],
can be represented by a Grassmann integral involving a Gaussian term
and a particular quartic term.
Although this representation has not yet led to any new rigorous results
concerning $F_G(\bw)$, it has led to important non-rigorous
insights into the behavior of $F_G(\bw)$ for large subgraphs
of a regular two-dimensional lattice
\cite{Caracciolo_matrixtree,JSS_forests}
--- insights that may yet be translatable into theorems
by exploiting the rigorous renormalization-group methods
developed in recent decades by mathematical physicists
(see e.g.\ \cite{Rivasseau_91}).

As in Section~\ref{subsec.matrixtree},
let $G=(V,E)$ be a finite undirected graph,
and let $L = L_G(\bw)$ be the Laplacian matrix for $G$
with edge weights $\bw = \{w_e\}_{e \in E}$.
Let us now introduce, at each vertex $i \in V$,
a pair of Grassmann variables $\psi_i$, $\psibar_i$.
All of these variables are nilpotent of order 2
(i.e., $\psi_i^2 = \psibar_i^2 = 0$),
anticommute, and obey the usual rules for Grassmann integration
\cite{Abdesselam_03,Zinn-Justin}.
Writing $\scrd(\psi,\psibar) = \prod_{i \in V} d\psi_i \, d\psibar_i$,
we have, for any matrix $A$,
\begin{equation}
   \int \scrd(\psi,\psibar) \; e^{\psibar A \psi}
   \;=\;
   \det A
\end{equation}
and more generally
\begin{eqnarray}
   & &  \!\!\!\!\!\!\!\!
        \int \! \scrd(\psi,\psibar) \; \psibar_{i_1} \psi_{j_1} \,\cdots\,
               \psibar_{i_r} \psi_{j_r} \, e^{\psibar A \psi}
     \nonumber \\
   & &  =\,
   \epsilon(i_1,\ldots,i_r|j_1,\ldots,j_r) \,
   \det A_{\setminus \{i_1,\ldots,i_r\},\{j_1,\ldots,j_r\}}
   \;,
\end{eqnarray}
where the sign $\epsilon(i_1,\ldots,i_r|j_1,\ldots,j_r) = \pm 1$
depends on how the vertices are ordered
but is always $+1$ when $(i_1,\ldots,i_r) = (j_1,\ldots,j_r)$.
These formulae allow us to rewrite the matrix-tree theorems in Grassmann form;
for instance, \reff{eq.principal} becomes
\begin{equation}
   \int \! \scrd(\psi,\psibar)
               \left( \prod_{\alpha=1}^r \psibar_{i_\alpha} \psi_{i_\alpha}
               \!\right)
               \, e^{\psibar L \psi}
   \;=
   \sum_{F \in {\cal F}(i_1,\ldots,i_r)} \, \prod_{e \in F} w_e   \;.
 \label{eq.principal.2}
\end{equation}

Let us next introduce,
for each connected (not necessarily spanning) subgraph
$\Gamma = (V_\Gamma, E_\Gamma)$ of $G$,
the object
\begin{equation}
   Q_\Gamma  \;=\;
   \left( \prod_{e \in E_\Gamma} w_e \right)
   \left( \prod_{i \in V_\Gamma} \psibar_i \psi_i \right)
   \;.
\end{equation}
(Note that each $Q_\Gamma$ is even and hence commutes with
the entire Grassmann algebra.)
Now consider an unordered family $\bGamma = \{\Gamma_1,\ldots,\Gamma_l\}$
with $l \ge 0$,
and let us try to evaluate an expression of the form
\begin{equation}
   \int \scrd(\psi,\psibar)
               \; Q_{\Gamma_1} \cdots Q_{\Gamma_l}
               \; e^{\psibar L \psi}
   \;.
\end{equation}
If the subgraphs $\Gamma_1,\ldots,\Gamma_l$
have one or more vertices in common,
then this integral vanishes on account of the nilpotency of the
Grassmann variables.
If, by contrast, the $\Gamma_1,\ldots,\Gamma_l$
are vertex-disjoint,
then \reff{eq.principal.2} expresses
$\int \! \scrd(\psi,\psibar) 
 \left( \prod_{k=1}^l \prod_{i \in V_{\Gamma_k}} \psibar_i \psi_i \right)
 e^{\psibar L \psi}$
as a sum over forests rooted at the vertices of
$V_{\bGamma} = \bigcup_{k=1}^l V_{\Gamma_k}$.
In particular, all the edges of
$E_{\bGamma} = \bigcup_{k=1}^l E_{\Gamma_k}$
must be absent from these forests,
since otherwise two or more of the root vertices would lie in the
same component (or one of the root vertices would be connected to itself
by a loop edge).
On the other hand, by adjoining the edges of $E_{\bGamma}$,
these forests can be put into one-to-one correspondence
with what we shall call \defword{$\bGamma$-forests}\/,
namely, spanning subgraphs $H$ in $G$ whose edge set contains $E_{\bGamma}$
and which, after deletion of the edges in $E_{\bGamma}$,
leaves a forest in which each tree component contains exactly one
vertex from $V_{\bGamma}$.
(Equivalently, a $\bGamma$-forest is a subgraph $H$
 with $l$ connected components in which each component
 contains exactly one $\Gamma_i$, and which does not contain any cycles
 other than those lying entirely within the $\Gamma_i$.
 Note, in particular, that a $\bGamma$-forest is a forest
 if and only if all the $\Gamma_i$ are trees.)
Furthermore, adjoining the edges of $E_{\bGamma}$
provides precisely the factor $\prod_{e \in E_\bGamma} w_e$.
Therefore
\begin{equation}
   \int \! \scrd(\psi,\psibar)
               \; Q_{\Gamma_1} \cdots Q_{\Gamma_l}
               \, e^{\psibar L \psi}
   \;=\;
   \sum_{H \in {\cal F}_\bGamma} \, \prod_{e \in H} w_e
 \label{eq.gamma}
\end{equation}
where the sum runs over all $\bGamma$-forests $H$.

We can now combine all the formulae \reff{eq.gamma} into a single
generating function, by introducing a variable $t_\Gamma$
for each connected subgraph $\Gamma$ of $G$.
Since $1 + t_\Gamma Q_\Gamma = e^{t_\Gamma Q_\Gamma}$, we have
\begin{equation}
   \int \! \scrd(\psi,\psibar)
            \; e^{\psibar L \psi + \sum\limits_\Gamma t_\Gamma Q_\Gamma}
\; = \!\!
       \sum_{\begin{scarray}
                \hbox{\scriptsize $\bGamma$ vertex-} \\
                \hbox{\scriptsize disjoint}
             \end{scarray}}
\!\!\!
       \Big( \prod_{\Gamma \in \bGamma} t_\Gamma \Big)
       \sum_{H \in {\cal F}_\bGamma} \, \prod_{e \in H} w_e
   \,.
 \label{eq.genfun1}
\end{equation}
We can express this in another way by interchanging
the summations over $\bGamma$ and $H$.
Consider an arbitrary spanning subgraph $H$
with connected components $H_1,\ldots,H_l$;
let us say that $\Gamma$ {\em marks}\/ $H_i$ (denoted $\Gamma \prec H_i$)
in case $H_i$ contains $\Gamma$ and contains no cycles other than those
lying entirely within $\Gamma$.  Define the weight
\begin{equation}
   W(H_i)  \;=\; \sum_{\Gamma \prec H_i} t_\Gamma   \;.
\end{equation}
Then saying that $H$ is a $\bGamma$-forest
is equivalent to saying that each of its components
is marked by exactly one of the $\Gamma_i$;
summing over the possible families $\bGamma$, we obtain
\begin{equation}
        \int \! \scrd(\psi,\psibar)
            \; e^{\psibar L \psi + \sum\limits_\Gamma t_\Gamma Q_\Gamma}
   \;\,= \!\!
       \sum_{\begin{scarray}
                H \, \hbox{\scriptsize spanning} \subseteq G \\
                H = (H_1,\ldots,H_l)
             \end{scarray}}
       \!\!
       \left( \prod_{i=1}^l W(H_i) \!\right)
       \prod_{e \in H} w_e
       \,.
 \label{eq.genfun2}
\end{equation}
This is our general combinatorial formula.
Extensions allowing prefactors
$\psibar_{i_1} \psi_{j_1} \cdots \psibar_{i_r} \psi_{j_r}$
are also easily derived.

Now consider the special case in which $t_\Gamma = t$ whenever
$\Gamma$ consists of a single vertex with no edges,
$t_\Gamma = u$ whenever $\Gamma$ consists of two vertices
linked by a single edge,
and $t_\Gamma = 0$ otherwise.
We have
\begin{multline}
     \int \! \scrd(\psi,\psibar)
        \, \exp\!\Big[
          \psibar L \psi \,
         +\, t \sum\limits_i \psibar_i \psi_i\,
         +\, u \sum\limits_{\< ij \>} w_{ij} \psibar_i \psi_i \psibar_j \psi_j
               \Big]
\\
   = \!\!\!\!\!\!\!\!
       \sum_{\begin{scarray}
                F \in {\cal F} \\
                F = (F_1,\ldots,F_l)
             \end{scarray}}
       \!\!\!\!\!\!
       \left( \prod_{i=1}^l \,  (t|V_{F_i}| + u|E_{F_i}|) \!\right)
       \,
       \prod_{e \in F} w_e
 \label{eq.fourfermion}
\end{multline}
where the sum runs over spanning forests $F$ in $G$
with components $F_1,\ldots,F_l$;
here $|V_{F_i}|$ and $|E_{F_i}|$ are, respectively,
the numbers of vertices and edges in the tree $F_i$.
[The four-fermion term
 $u \sum_{\< ij \>} w_{ij} \psibar_i \psi_i \psibar_j \psi_j$
 can equivalently be written, using nilpotency of the Grassmann variables,
 as $-(u/2) \sum_{i,j} \psibar_i \psi_i L_{ij} \psibar_j \psi_j$.]
If $u=0$, this formula represents vertex-weighted spanning forests
as a determinant (``massive fermionic free field'')
\cite{Duplantier_88,Biggs_93}.
More interestingly,
since $|V_{F_i}| - |E_{F_i}| = 1$ for each tree $F_i$,
we can take $u=-t$
and obtain the generating function of {\em unrooted}\/ spanning forests
with a weight $t$ for each component.
This is furthermore equivalent to giving each edge $e$ a weight $w_e/t$,
and then multiplying by an overall prefactor $t^{|V|}$.
We have therefore proven:

\begin{prop}
  \label{prop.forests}
Let $G=(V,E)$ be a finite undirected graph,
let $L$ be the Laplacian matrix for $G$
with edge weights $\bw = \{w_e\}_{e \in E}$,
and let $F_G$ be the generating polynomial of spanning forests in $G$.
Then
\begin{equation}
     \int \! \scrd(\psi,\psibar)
        \, \exp\!\Big[
          \psibar L \psi \,
         +\, t \sum\limits_i \psibar_i \psi_i\,
         +\, {t \over 2}
           \sum\limits_{i,j} \psibar_i \psi_i L_{ij} \psibar_j \psi_j
               \Big]
     \;=\;
     t^{|V|} F_G(\bw/t)
     \;.
\end{equation}
\end{prop}

\noindent
This representation of unrooted spanning forests is the translation
to generating functions and Grassmann variables of a little-known
but important paper by Liu and Chow \cite{Liu_81}.

It is an open question whether there exists some analogue of
Proposition~\ref{prop.forests} for the object dual to $F_G(\bw)$,
namely the generating polynomial $C_G(\bv)$ of connected spanning subgraphs,
when the graph $G$ is non-planar.

It would also be interesting to know whether the foregoing identities
are in any way related to the forest-root formula of Brydges and Imbrie
\cite{Brydges-Imbrie_03a,Brydges-Imbrie_03b,Imbrie_03,Imbrie_04}.

\section{Absence of zeros at large \mbox{\protect\boldmath $|q|$}}
   \label{sec.largeq}

Combinatorialists have long been interested in the real or complex zeros
of the chromatic polynomial $P_G(q)$.
A fair number of interesting theorems have by now been proven,
notably concerning the real zeros,
but even here a vast number of open problems remain;
moreover, very little is known rigorously concerning the complex zeros.
(See \cite{Jackson_03} for an excellent recent survey
 treating both real and complex zeros.)

In this section I~would like to discuss one aspect of the complex-zero
problem where some progress has recently been made by exploiting
methods from mathematical statistical mechanics:
namely, the absence of zeros at large enough $|q|$.

\subsection{Bounds in terms of maximum degree and its relatives}

The bounds I~want to discuss \cite{Sokal_chromatic_bounds}
apply, in fact, not only to the chromatic polynomial
but to the multivariate Tutte polynomial $Z_G(q,\bv)$ throughout
the ``complex antiferromagnetic regime'' $|1+v_e| \le 1$.
These bounds come in several variants, but here is a typical one:

\begin{theorem}
  \label{thm.sokal}
Let $G=(V,E)$ be a loopless finite undirected graph
of maximum degree $\Delta$.
Suppose that the edge weights $\bv = \{ v_e \}_{e \in E}$
satisfy $|1 + v_e| \le 1$ for all $e$.
Let $v_{\rm max} = \max\limits_{e \in E} |v_e|$.
Then all the zeros of $Z_G(q, \bv)$
lie in the disc $|q| < 7.963907 v_{\rm max} \Delta$.
\end{theorem}

\noindent
For the chromatic polynomial ($v_e = -1$ for all $e$)
one obtains the immediate corollary:

\begin{cor}
   \label{cor.sokal}
Let $G=(V,E)$ be a loopless finite undirected graph
of maximum degree $\Delta$.
Then all the zeros of the chromatic polynomial $P_G(q)$
lie in the disc $|q| < 7.963907 \Delta$.
\end{cor}

\noindent
Corollary~\ref{cor.sokal} answers in the affirmative a question posed by
Brenti, Royle and Wagner \cite[Question 6.1]{Brenti_94},
generalizing an earlier conjecture of
Biggs, Damerell and Sands \cite{Biggs_72}
limited to regular graphs.

Of course, the constant 7.963907 is an artifact of the proof;
it is presumably far from sharp.
The linear dependence on $\Delta$ is, however, best possible,
since the complete graph $K_{\Delta+1}$
has chromatic roots $0,1,2,\ldots,\Delta$.
Furthermore (and perhaps surprisingly),
the complete graph $K_{\Delta+1}$ is {\em not}\/
the extremal graph for this problem,
and a bound $|q| \le \Delta$ is {\em not}\/ possible.
In fact, a non-rigorous (but probably rigorizable) asymptotic analysis,
confirmed by numerical calculations, shows \cite{Salas-Sokal_bipartite}
that the complete bipartite graph $K_{\Delta,\Delta}$
has a chromatic root $\alpha \Delta + o(\Delta)$, where
$\alpha = - 2 / W(-2/e)  \approx 0.678345 + 1.447937 i$;
here $W$ denotes the principal branch of the Lambert $W$ function
(the inverse function of $w \mapsto w e^w$) \cite{Corless_96}.
So the constant in Corollary~\ref{cor.sokal} cannot be better than
$|\alpha| \approx 1.598960$.

A complete proof of Theorem~\ref{thm.sokal} and its variants
can be found in \cite{Sokal_chromatic_bounds}
(see also \cite{Borgs_lectures});
here I~would like simply to sketch the method,
which I~hope will be of wider use.

It is immediately apparent that Theorem~\ref{thm.sokal} and
Corollary~\ref{cor.sokal} are ``hard'' results in the sense of
Section~\ref{sec.complex}:
the zero-free region depends on the maximum degree $\Delta$.
The proof of these results follows the plan sketched there:
namely, the multivariate Tutte polynomial on the graph $G=(V,E)$
is first mapped onto a gas of nonoverlapping ``polymers'' on $V$,
i.e.\ a hard-core lattice gas on the intersection graph of $V$;
then this hard-core lattice gas is controlled by means of the
theorems mentioned at the end of Section~\ref{sec.complex}.

The mapping to a polymer gas is really quite simple;
indeed, it is nearly trivial.
The definition \reff{eq.def.Ztilde} writes
the multivariate Tutte polynomial $\Ztilde_G(q,\bv)$
as a sum over spanning subgraphs $(V,A)$.
Let us perform this sum in two stages:
First we sum over partitions of $V$ into disjoint nonempty subsets
$V_1,\ldots,V_k$ ($k \ge 1$ arbitrary);
these subsets will correspond to the vertex sets of the connected components
of $(V,A)$.
Then we sum over the ways of choosing edges within each component
so as to make it connected.
{}From \reff{eq.def.Ztilde} we see that a component $V_i = S$
will get a weight
\begin{subeqnarray}
   w(S)   & = &
   q^{1-|S|}
   \!\!\!\!   \sum\limits_{\begin{scarray}
                               B \subseteq E  \\
                               (S,B) \,\text{\scriptsize connected}
                           \end{scarray}}
               \prod\limits_{e \in B} v_e
 \slabel{eq.weights.polymer.a}
       \\
   & = &  q^{1-|S|} \, C_{G|S}(\bv)
 \label{eq.weights.polymer}
\end{subeqnarray}
where $G|S$ denotes the induced subgraph.
Note, in particular, that (provided $G$ is loopless)
any set $S$ of cardinality 1 gets weight $w(S) = 1$.
So we need not consider such sets explicitly:
it suffices to define the ``polymers''
to be the sets $V_i$ of cardinality 2 or more;
we then know that any vertex not covered by a polymer
must be covered by a set $V_i$ of cardinality 1.
The polymers are therefore an arbitrary family (possibly empty)
of disjoint sets $S_i \subseteq V$ of cardinality 2 or more.
We have thus written $\Ztilde_G(q,\bv)$
as an independent-set polynomial (= hard-core lattice gas)
of the form \reff{eq.Zhardcore}
for a graph $\widehat{G} = (\widehat{V},\widehat{E})$
whose vertices are the subsets of $V$ of cardinality 2 or more,
and whose edges are the pairs with nonempty intersection.
The fugacities are given by \reff{eq.weights.polymer}.

The usefulness of this representation for proving the absence of zeros
at large $|q|$ comes from the fact that the fugacities $w(S)$
are all suppressed by powers of $q^{-1}$ (since $|S| \ge 2$),
hence are small for large $|q|$.
Furthermore, this suppression operates more strongly for larger polymers.
This raises the hope that, when $|q|$ is sufficiently large,
the model lies in a ``low-fugacity'' regime
where the methods mentioned in Section~\ref{sec.complex}
can be brought to bear.  This is in fact the case.

Of course, I~have sloughed over one crucial point:
the sum over $B$ in \reff{eq.weights.polymer.a}.
It is here, not surprisingly, that the maximum degree $\Delta$ enters
(along with $v_{\rm max}$).
What one needs, it turns out, is a bound on the number of
connected subgraphs $H \subseteq G$ containing a specified vertex $x$
and having a specified number of vertices $m=|S|$.
What is easier to prove, however, is a bound on the number of
connected subgraphs $H \subseteq G$ containing a specified vertex $x$
and having a specified number of {\em edges}\/
\cite[Section 4.2]{Sokal_chromatic_bounds}
\cite{Jackson-Sokal_maxflow}.
The relation between the two is far from obvious,
but it turns out that in the complex antiferromagnetic region
$|1 + v_e| \le 1$ (and only there!),
a beautiful inequality due to Oliver Penrose \cite{Penrose_67}
allows connected subgraphs to be bounded by {\em trees}\/.\footnote{
   See \cite[Section~4.1]{Sokal_chromatic_bounds}
   and \cite[Section~2.2]{Scott-Sokal}
   for further discussion.
}
The upshot is that the sum over $B$ in \reff{eq.weights.polymer.a}
can be controlled, with the result that $w(S)$
satisfies a uniform bound that is exponentially decaying in $|S|$,
provided that $|q|$ is large enough compared to $\Delta$.
This suffices to ensure the convergence of the ``polymer expansion'',
and hence the nonvanishing of $\Ztilde_G(q,\bv)$.

I~believe it should be possible to extend this result
beyond the complex antiferromagnetic regime,
by bounding directly the number of connected subgraphs
having a specified number of vertices \cite{Jackson-Sokal_vertices}.
But the resulting bound on $|q|$ will no longer be
linear in $\Delta$ and $v_{\rm max}$;
rather, it will behave roughly like $v_{\rm max}^{\Delta/2}$.
And this is not simply an artifact of the method of proof:
in the $q$-state Potts {\em ferromagnet}\/ ($v > 0$, $q > 0$)
on the the simple hypercubic lattice $\Z^d$ with nearest-neighbor edges
(hence $\Delta=2d$),
there is a first-order phase transition at
$v_t(q)   =  q^{1/d} \, [1 + O(1/q)]$
for all sufficiently large $q$
\cite{Martirosian_86,Laanait_86,Kotecky_90,Laanait_91,Borgs_91}.
As discussed in Section~\ref{sec.complex},
this means that for $v$ near $v_t(q)$
there will be complex zeros of $Z_G(q,v)$
for large subgraphs $G \subset \Z^d$.

By an {\em ad~hoc}\/ (and aesthetically unsatisfying) method,
Theorem~\ref{thm.sokal} and Corollary~\ref{cor.sokal}
can be strengthened so as to bound the roots,
not in terms of the maximum degree $\Delta$,
but in terms of the {\em second-largest}\/ degree $\Delta_2$
\cite[Section~6]{Sokal_chromatic_bounds}.
For simplicity let me state only the result for chromatic polynomials:

\begin{cor}
   \label{cor.sokal2}
Let $G=(V,E)$ be a loopless finite undirected graph
of second-largest degree $\Delta_2$.
Then all the zeros of the chromatic polynomial $P_G(q)$
lie in the union of the discs
$|q| < 7.963907 \Delta_2$ and $|q-1| < 7.963907 \Delta_2$.
In particular, they all lie in the disc $|q| < 7.963907 \Delta_2 + 1$.
\end{cor}

\noindent
It should be stressed that ``second-largest''  in Corollary \ref{cor.sokal2}
{\em cannot}\/ be replaced by ``third-largest''.
Indeed, the generalized theta graphs $\Theta^{(s,p)}$
\cite{gen_theta,Sokal_chromatic_roots} constitute a family
of planar (in fact, series-parallel) graphs
in which all but two vertices have degree 2,
yet their chromatic roots are dense in the region
$\{ q \in \C \colon\;  |q-1| \ge 1 \}$ \cite{Sokal_chromatic_roots}.\footnote{
   The graph $\Theta^{(s,p)}$ consists of a pair of endvertices
   connected by $p$ internally disjoint paths each of length $s$.
}
So, with {\em one}\/ large-degree vertex, the chromatic roots remain bounded;
but with {\em two}\/, all hell can break loose.

\subsection{Bounds in terms of maxmaxflow?}

Alas, all of the foregoing bounds have a severe defect:
they are {\em unnatural}\/,
because the multivariate Tutte polynomial $Z_G(q,\bv)$
factorizes over blocks [cf.\ \reff{eq.blocks}],
while $\Delta$ and $\Delta_2$ can grow arbitrarily large
when blocks are glued together at a cut vertex.
As a consequence, the bounds obtained from Theorem~\ref{thm.sokal}
and its variants are sometimes very far from sharp:
for instance, a tree can have arbitrarily large $\Delta$ and $\Delta_2$,
but its chromatic roots are only 0 and 1!
Of course, this defect can be cured by the {\em ad~hoc}\/ technique
of applying Theorem~\ref{thm.sokal} ff.\ 
to each block of $G$ rather than directly to $G$.
But the deeper significance of this remark
is that the quantities $\Delta$ and $\Delta_2$
are {\em unnatural}\/ for studying the multivariate Tutte polynomial
because they are {\em not matroidal}\/.

One would like, therefore, to find a graph invariant that
is smaller than $\Delta$ and $\Delta_2$,
that ``trivializes over blocks'' (and ideally generalizes to matroids),
and that is strong enough to bound the roots of $P_G(q)$ and $Z_G(q,\bv)$.
I~have a candidate:  it is the {\em maxmaxflow}\/
\cite{Shrock_99a,Sokal_chromatic_bounds,Jackson-Sokal_maxflow},
defined as follows:
For $x,y \in V$ with $x \neq y$, let $\lambda_G(x,y)$ be the max flow
from $x$ to $y$:
\begin{subeqnarray}
   \lambda_G(x,y)
      & = &  \hbox{max \# of edge-disjoint paths from $x$ to $y$} \\
      & = &  \hbox{min \# of edges separating $x$ from $y$}
\end{subeqnarray}
The {\em maxmaxflow}\/ $\Lambda(G)$ is then defined by
\begin{equation}
   \Lambda(G)   \;=\;
   \max\limits_{\begin{scarray}
                   x,y \in V \\
                   x \neq y
                \end{scarray}}
   \lambda_G(x,y)   \;.
\end{equation}
[Note the contrast with the edge-connectivity, which is the {\em minimum}\/
of $\lambda_G(x,y)$ over $x \neq y$.]
Clearly $\lambda_G(x,y) \le \min[d_G(x), d_G(y)]$,
so that
\begin{equation}
   \Lambda(G) \;\le\; \Delta_2(G)  \;.
\end{equation}
Furthermore, it is easy to see that maxmaxflow ``trivializes over blocks''
in the sense that $\Lambda(G) = \max_{1 \le i \le b} \Lambda(G_i)$
where $G_1,\ldots,G_b$ are the blocks of $G$.

An apparently very different quantity can be defined via cocycle bases.
For $X,Y$ disjoint subsets of $V$,
let $E(X,Y)$ denote the set of edges in $G$ between $X$ and $Y$.
A \defword{cocycle}\/ of $G$ is a set $E(X,X^c)$
where $X \subseteq V$ and $X^c \equiv V \setminus X$.
It is well-known that the cocycles of $G$
form a vector space over $GF(2)$ with respect to symmetric difference;
this is called the \defword{cocycle space}\/ of $G$.
Let $\widetilde{\Lambda}(G)$ be the minmax cardinality of the
cocycles in a basis, i.e.
\begin{equation}
   \widetilde{\Lambda}(G)  \;=\;
   \min\limits_{\scrb}  \max\limits_{C\in \scrb} |C|
   \;,
\end{equation}
where the min runs over all bases $\scrb$ of the cocycle space of $G$,
and the max runs over all the cocycles $C$ in the basis $\scrb$.
Since one special class of cocycle bases consists
of taking the stars
$C(x) = E(\{x\}, \{x\}^c)$
for all but one of the vertices in each component of $G$, we clearly have
\begin{equation}
   \widetilde{\Lambda}(G) \;\le\; \Delta_2(G)  \;.
\end{equation}

The relationship, if any, between maxmaxflow and cocycle bases
is perhaps not obvious at first sight.
But Bill Jackson and I have proven \cite{Jackson-Sokal_maxflow} that
\begin{equation}
   \Lambda(G)  \;=\;  \widetilde{\Lambda}(G) \;.
\end{equation}
The two definitions thus give dual approaches to the same quantity.

All this theory, extends in fact, to the more general situation of a
finite undirected graph $G$ equipped with
nonnegative real edge weights $\bw = \{ w_e \} _{e \in E}$.
(In the application to the multivariate Tutte polynomial,
 one might want to take $w_e = |v_e|$ or something similar.)

The ultimate goal of this work is to prove (or disprove!)\ 
strengthenings of Theorem~\ref{thm.sokal} ff.\ 
in which maximum degree or second-largest degree are replaced by maxmaxflow.
For instance, for chromatic polynomials one would like to prove the following:

\begin{conj}
   \label{conj.chrom.1}
There exist universal constants $C(\Lambda) < \infty$
such that all the chromatic roots (real or complex)
of all loopless graphs of maxmaxflow $\Lambda$
lie in the disc $|q| \le C(\Lambda)$.
Indeed, I conjecture that $C(\Lambda)$ can be taken to be linear in $\Lambda$.
\end{conj}

\noindent
It is natural to try to prove Conjecture~\ref{conj.chrom.1}
by modifying the arguments of \cite{Sokal_chromatic_bounds}
so as to decompose a spanning subgraph of $G$ into
its {\em 2-connected}\/ (rather than merely {\em connected}\/) components.
Such an argument will require, at a minimum,
a good bound on the number of 2-connected subgraphs $H \subseteq G$
containing a specified edge $e$ and having a specified number of
vertices or edges.
Bill Jackson and I have found some bounds of this type
\cite{Jackson-Sokal_maxflow}.
But this is only a first step,
and the other obstacles in generalizing this proof ---
such as controlling the interaction between the 2-connected components ---
have yet to be overcome.

It is worth mentioning that $\widetilde{\Lambda}$
can be defined also for binary matroids,
generalizing the definition for graphs.
But we do not yet know whether it is strong enough to give a good
(i.e., singly exponential) bound on the number of 2-connected
deletion minors containing a specified element $e$
and having a specified size.
Nor do we know whether it is useful in bounding the zeros
of the matroid chromatic polynomial $\Ztilde_M(q,-1)$.

\subsection{Some further questions}

There is no compelling reason, {\em a priori}\/,
to limit attention in the complex $q$-plane to discs centered at the origin;
many other types of regions are of interest.
For instance, the following conjecture might conceivably be true:
\begin{conj}
   \label{conj.chrom.2}
All the chromatic roots (real or complex)
of all loopless graphs of maxmaxflow $\Lambda$
lie in the half-plane $\real q \le \Lambda$.
\end{conj}
By factoring $P_G(q)/q$ and using the fact that its roots are either real
or occur in complex-conjugate pairs, one easily deduces that
the truth of Conjecture ~\ref{conj.chrom.2} would imply:
\begin{conj}
   \label{conj.chrom.3}
If $G$ has maxmaxflow $\Lambda$,
then $P_G(q)/q$ and all its derivatives are nonnegative at $q=\Lambda$.
[The same then holds also for $P_G(q)$ and all its derivatives.]
\end{conj}
Since the truth of Conjecture ~\ref{conj.chrom.3} would imply
that $P_G(q)$ has all nonnegative Taylor coefficients at $q=\Lambda$,
and since $P_G(q)$ is not identically zero (when $G$ is loopless),
Conjecture ~\ref{conj.chrom.3} would imply:
\begin{conj}
   \label{conj.chrom.4}
If $G$ is a loopless graph of maxmaxflow $\Lambda$,
then $P_G(q) > 0$ for $q > \Lambda$.
In particular, $P_G(q)$ has no roots in the real interval $(\Lambda,\infty)$.
\end{conj}

Even the weaker versions of these conjectures with $\Lambda$
replaced by $\Delta_2$ or $\Delta$ are open!
So one has the $3 \times 3$ matrix of conjectures
\begin{eqnarray*}
\begin{array}{ccccc}
    \ref{conj.chrom.2} & \Longrightarrow & \ref{conj.chrom.2}' & \Longrightarrow
 &
                                      \ref{conj.chrom.2}''                    \\
    \Downarrow    &       &  \Downarrow    &       &  \Downarrow         \\
    \ref{conj.chrom.3} & \Longrightarrow & \ref{conj.chrom.3}' & \Longrightarrow
 &
                                      \ref{conj.chrom.3}''                    \\
    \Downarrow    &       &  \Downarrow    &       &  \Downarrow         \\
    \ref{conj.chrom.4} & \Longrightarrow & \ref{conj.chrom.4}' & \Longrightarrow
 &
                                      \ref{conj.chrom.4}''
\end{array}
\end{eqnarray*}
in which the strongest might well be true
and the weakest is still an open question.\footnote{
   {\bf STOP PRESS!!!}
   Gordon Royle (private communication)
   has found a family of counterexamples to
   Conjectures~\ref{conj.chrom.2} and \ref{conj.chrom.2}${}'$.
   His graphs are planar (but not series-parallel)
   and all vertices but one have degree 3 (hence $\Lambda = \Delta_2 = 3$).
   Included in his family are a 47-vertex graph with chromatic roots
   $q \approx 3.0129950712 \pm 0.8089628639\,i$,
   a 95-vertex graph with chromatic roots
   $q \approx 3.0536525915 \pm 0.7547530551\,i$,
   a 191-vertex graph with chromatic roots
   $q \approx 3.07174237056 \pm 0.7105232675\,i$,
   and a 383-vertex graph with chromatic roots
   $q \approx 3.0766232972 \pm 0.6746120243\,i$.
}
It might be useful to seek counterexamples,
either by systematic calculation on small graphs up to $\approx 20$ vertices
or, perhaps more fruitfully, by constructing suitable infinite families
of graphs.
It goes without saying that I~have no idea how to prove
any of these conjectures, or even any compelling reason to believe
they are true (other than my inability thus far to find counterexamples).

We can pose these questions most generally as follows:
Let $\scrg$ be a class of finite graphs,
and let $\scrv$ be a subset of the complex plane.
Then we can ask about the sets
\begin{eqnarray}
   S_1(\scrg, \scrv)   & = &
      \bigcup_{G \in \scrg} \;
      \bigcup_{v \in \scrv} \;
      \{q \in \C \colon\;  Z_G(q,v) \,=\, 0 \}               \\[2mm]
   S_2(\scrg, \scrv)   & = &
      \bigcup_{G \in \scrg} \;
      \bigcup_{\bv \colon v_e \in \scrv \; \forall e}
      \{q \in \C \colon\;  Z_G(q, \bv) \,=\, 0 \}
\end{eqnarray}
Among the interesting cases are
the chromatic polynomials $\scrv = \{-1\}$,
the antiferromagnetic regime $\scrv = [-1,0]$,
and the complex antiferromagnetic regime
$\scrv = \scra \equiv \{v \in \C \colon\; |1 + v| \le 1 \}$.

The theorems discussed in this section show that
the sets $S_2(\scrg_r, \scra)$ and $S_2(\scrg'_r, \scra)$ are {\em bounded}\/,
where $\scrg_r$ (resp.\ $\scrg'_r$) 
is the set of all graphs whose maximum degree (resp.\ second-largest degree)
is less than or equal to $r$.
But boundedness is a rather weak statement about a set;
one would like to learn more about its location in the complex plane.

In the other direction, I~have recently shown \cite{Sokal_chromatic_roots}
that if $\scrg$ is the family of all generalized theta graphs $\Theta^{(s,p)}$
\cite{gen_theta,Sokal_chromatic_roots},
then $S_1(\scrg, \{-1\})$ is dense in the region
$\{ q \in \C \colon\;  |q-1| \ge 1 \}$.
Moreover, if $\scrg'$ denotes the family of
joins of $\Theta^{(s,p)}$ with the complete graph $K_2$,
then $S_1(\scrg \cup \scrg', \{-1\})$
is dense in the entire complex plane.\footnote{
   The {\em join}\/ of $G$ and $H$, denoted $G+H$,
   is the graph obtained from the disjoint union of $G$ and $H$
   by adding one edge connecting each pair of vertices
   $x \in V(G)$, $y \in V(H)$.
}
The idea of the construction in \cite{Sokal_chromatic_roots}
is to exploit the reduction formulae of Section~\ref{subsec.2-rooted}:
suitably chosen 2-rooted subgraphs (in this case single edges)
are concatenated in order to create a larger 2-rooted subgraph
that has some desired value of the effective coupling $v_{\rm eff}$.
Such a construction is obviously inapplicable to 3-connected graphs,
in which 2-rooted subgraphs of more than a single edge cannot occur.
It is thus reasonable to ask:
What is the closure of the set of chromatic roots of
{\em 3-connected}\/ graphs?
The answer, alas, is still the whole complex plane \cite{Jackson_private}.
To see this, consider the graphs $\Theta^{(s,p)} + K_n$
for fixed $n$ and varying $s,p$:
they are $(n+2)$-connected,
and their chromatic roots (taken together)
are dense in the region $|q - (n+1)| \ge 1$.
In particular, considering both $n$ and $n+2$,
the chromatic roots are dense in the whole complex plane.
So even arbitrarily high connectedness does not, by itself,
stop the chromatic roots from being dense in the whole complex plane.

On the other hand, the graphs $\Theta^{(s,p)} + K_n$ are non-planar
whenever $n \ge 1$ and $p \ge 3$.
This suggests posing a more restricted question:
\begin{qn}
What is the closure of the set of chromatic roots for 3-connected
 {\em planar}\/ graphs?
\end{qn}
Here the answer may well be much smaller than
$\C$ or $\C \setminus \{|q-1| < 1\}$.
But it will not be a bounded set \cite{Jackson_private},
since the bipyramids $C_n + \bar{K}_2$ are 4-connected plane triangulations,
but their chromatic roots are unbounded \cite{Read_88,Shrock_97a,Shrock_97c}.

\subsection{A final remark}

The ultimate goal of this research is to prove theorems of the type
``such-and-such graph property (or conjunction of properties)
implies such-and-such bound on the chromatic roots''.
But there is a virtually unlimited number of possible assertions
of this type;  the hard part is to figure out which ones are true!

There is thus, in my opinion, much room for numerical experiment
as a means for obtaining intuition about which graph properties
(e.g.\ maximum degree and its relatives, girth, planarity,
 connectivity, \ldots\!\!\!)
affect which aspects of the chromatic roots
(e.g.\ large positive real roots, large positive or negative real part,
 large imaginary part, etc.).
Some information can be obtained by systematic calculation on all graphs
of a specified type up to as many vertices as can be handled
--- or, for suitable classes (e.g.\ cubic graphs),
by random generation of typical graphs
with a rather larger number of vertices ---
in an attempt to correlate properties of the graph with
properties of the chromatic roots.
Guided by this information, one may then choose to focus on special classes
of graphs that potentially exhibit some type of ``extremal'' behavior,
and either pursue the numerical calculation to larger graphs
or, better yet, find infinite families of graphs
that exhibit such ``extremal'' behavior
and for which one can compute exactly the asymptotic behavior
of the chromatic roots.

It is important to keep in mind that many asymptotic properties
of chromatic roots are achieved very slowly as the number of vertices grows:
for instance, the roots of cubic graphs move quite slowly
towards negative real part \cite{Read_91},
and the roots of the generalized theta graphs $\Theta^{(s,p)}$
fill out the complex plane extremely slowly as $s,p \to \infty$
\cite{gen_theta}.
As a result, numerical experiments can be deceptive unless
interpreted with extreme care.
It is for this reason that exact computations with
well-chosen infinite families of graphs are particularly valuable.

\thankyou{
   Many of the results discussed here were obtained in joint work
   with Jason Brown, Sergio Caracciolo, Young-Bin Choe, Carl Hickman,
   Bill Jackson, Jesper Jacobsen, James Oxley, Gordon Royle, Jes\'us Salas,
   Hubert Saleur, Alex Scott, Andrea Sportiello and Dave Wagner.
   I~thank them all for the pleasure of doing mathematics together,
   as well as for permission to plagiarize here some of our
   collective prose.
   I~also thank Abdelmalek Abdesselam, Norman Biggs, Graham Brightwell,
   David Brydges, Geoffrey Grimmett, John Imbrie, Marco Polin, Robert Shrock,
   Jan van den Heuvel, Dominic Welsh and Geoff Whittle
   for many helpful discussions;
   and Bill Jackson, Mark Jerrum, Gordon Royle, Dave Wagner and Geoff Whittle
   for helpful comments on a first draft of this paper.
   
   The author's research was supported in part
   by U.S.\ National Science Foundation grants
   PHY--0099393 and PHY--0424082
   and by U.K.\ Engineering and Physical Sciences Research Council
   grant GR/S26323/01.
}

%
%
%
%

\myaddress
\end{document}